\documentclass[a4paper]{article}

\RequirePackage[utf8]{inputenx}
\RequirePackage[T1]{fontenc}

\RequirePackage[final]{microtype}
\usepackage{textcomp}

\usepackage{amssymb}
\usepackage{mathtools}
\usepackage{amsthm}
\usepackage{thmtools}

\usepackage{graphicx}
\graphicspath{{figures/}}

\usepackage{array}
\usepackage{longtable}
\usepackage{booktabs}

\usepackage{sectsty}
\sectionfont{\raggedright}

\usepackage[obeyFinal]{todonotes}
\usepackage{xspace}

\usepackage{varioref}
\usepackage[pdfusetitle, colorlinks, allcolors = uiolink]{hyperref}
\usepackage[nameinlink, capitalize, noabbrev]{cleveref}

\usepackage[backend     = biber,
            sortcites   = true,
            giveninits  = true,
            maxbibnames = 5,
            doi         = false,
            isbn        = false,
            url         = false,
            sortlocale  = nb_NO,
            style       = alphabetic]{biblatex}
\DeclareNameAlias{sortname}{family-given}
\DeclareNameAlias{default}{family-given}
\DeclareFieldFormat[article]{volume}{\bibstring{jourvol}\addnbspace#1}
\DeclareFieldFormat[article]{number}{\bibstring{number}\addnbspace#1}
\renewbibmacro*{volume+number+eid}
{
    \printfield{volume}
    \setunit{\addcomma\space}
    \printfield{number}
    \setunit{\addcomma\space}
    \printfield{eid}
}

\declaretheorem[style = plain, numberwithin = section]{theorem}

\declaretheorem[style = plain,      sibling = theorem]{lemma}
\declaretheorem[style = plain,      sibling = theorem]{proposition}
\declaretheorem[style = plain,      sibling = theorem]{conjecture}

\declaretheorem[style = definition, sibling = theorem]{example}

\numberwithin{equation}{section}
\numberwithin{figure}{section}

\renewcommand{\Re}{\operatorname{Re}}
\renewcommand{\Im}{\operatorname{Im}}
\DeclareMathOperator{\rank}{rank}
\DeclareMathOperator{\corank}{corank}
\DeclareMathOperator{\Sing}{Sing}
\DeclareMathOperator{\Bl}{Bl}

\DeclarePairedDelimiter{\linspan}{\langle}{\rangle}
\DeclarePairedDelimiter{\abs}{\lvert}{\rvert}

\newcommand{\R}{\mathbb{R}}    
\newcommand{\C}{\mathbb{C}}    
\renewcommand{\P}{\mathbb{P}}  
\newcommand{\G}{\mathbb{G}}    

\newcommand{\x}{\mathbf{x}}
\newcommand{\y}{\mathbf{y}}

\newcommand{\e}{\mathbf{e}}
\newcommand{\ii}{\boldsymbol{i}}

\newcommand{\locus}[1]{rank-\(#1\) locus}
\newcommand{\point}[1]{rank-\(#1\) point}
\newcommand{\pointsec}[1]{rank-#1 point}
\newcommand{\node}[1]{rank-\(#1\) node}
\newcommand{\quadric}[1]{rank-\(#1\) quadric}
\newcommand{\matrices}[1]{rank-\(#1\) matrices}

\newcommand{\MXkern}{\ensuremath{\mkern-7mu}}
\newcommand{\Pnkern}{\ensuremath{\mkern-4mu}}
\newcommand{\Wkern}{\ensuremath{\mkern-4mu}}
\newcommand{\Hkern}{\ensuremath{\mkern-2mu}}

\newcommand{\dash}{\textthreequartersemdash\xspace}
\newcommand{\V}{\mathcal{V}\mkern-1mu}

\newcommand{\T}{\mathsf{T}}
\definecolor{uiolink}{HTML}{0B5A9D}

\title{Determinantal quartic surfaces\texorpdfstring{\\}{ }with a definite Hermitian representation}
\author{Martin Helsø}

\begin{document}

\maketitle

\begin{abstract}
    \noindent
    We give a bound on the number of isolated, essential singularities
    of determinantal quartic surfaces in 3-space.
    We also provide examples
    of different configurations of real singularities
    on quartic surfaces with a definite Hermitian determinantal representation,
    and conjecture an extension of a theorem by Degtyarev and Itenberg.
\end{abstract}

\section{Introduction}

Representing a polynomial as the determinant of a linear matrix
is a problem dating back at least to Hesse \cite{Hes44}.
Determinantal representations have applications in areas such as
linear algebra, operator theory, convex optimisation and algebraic geometry.

A homogeneous polynomial $f \in \C[x_0, \ldots, x_n]_d$ of degree~$d$
has a \emph{determinantal representation} if $f(\x) = \det(M_d(\x))$,
where
\begin{equation}
    \label{eq:hermitian-representation}
    M_d(\x)
    \coloneqq
    M_d(x_0, \ldots, x_n)
    \coloneqq
    M_{d, 0} x_0 + \cdots + M_{d, n} x_n
\end{equation}
for some $(d \times d)$-matrices~$M_{d, 0}, \ldots, M_{d, n}$.
We say that the representation is \emph{Hermitian}
if $M_{d, 0}, \ldots, M_{d, n}$ are Hermitian matrices,
and it is \emph{symmetric} if $M_{d, 0}, \ldots, M_{d, n}$ are symmetric matrices.
The hypersurface $\V(f) \subset \C\P^n$ is called \emph{determinantal}
if $f$ possesses a determinantal representation;
$V(f)$ is called a \emph{symmetroid} if the representation is symmetric.

Let $f \in \R[x_0, \ldots, x_n]_d$ be a real polynomial
with a Hermitian determinantal representation~\eqref{eq:hermitian-representation}.
The representation is \emph{definite} if the matrix~$M_d(\e)$ is positive definite
for some point~$\e \in \R\P^n$\Pnkern.
The eigenvalues of a Hermitian matrix are real.
It follows that every real line through $\e$
only meets the hypersurface $\V(f) \subset \C\P^n$ in real points.
A polynomial with this property is called \emph{hyperbolic with respect to $\e$}.
The connected component of $\e$ in $\R\P^n \setminus V_\R(f)$
is called the \emph{hyperbolicity cone} of $f$ with respect to $\e$.

If $M_{d, 0}, \ldots, M_{d, n}$ in \eqref{eq:hermitian-representation} are real, symmetric matrices,
then the set
\begin{equation*}
    \big\{
        \x \in \R\P^n
        \mid
        M_d(\x) \text{ is semidefinite}
    \big\}
\end{equation*}
is called a \emph{spectrahedron}.
It is easy to see that all spectrahedra are hyperbolicity cones.
The converse statement,
that all hyperbolicity cones are spectrahedra,
is called the generalised Lax conjecture
and is an object of much interest.
A partial result is the Helton--Vinnikov theorem,
which implies that all hyperbolicity cones in $\R\P^2$
are spectrahedra \cite{HV07, LPR05}.
Let $f \in \R[x_0, \ldots, x_n]_d$ be polynomial
with a definite Hermitian representation~\eqref{eq:hermitian-representation},
and suppose that $\e \in \R\P^n$ is such that $M_d(\e)$
is positive definite.
We note that the hyperbolicity cone of $f$ with respect to $\e$
is a spectrahedron, by following \cite[Corollary~5.1]{PV13}:
We can write $M_d(\x) = A_d(\x) + \ii B_d(\x)$,
where $A_d(\x)$ is real symmetric and $B_d(\x)$ is real antisymmetric.
We define the real symmetric $(2d \times 2d)$-matrix
\begin{equation}
    \label{eq:associated-sym}
    A_{2d}(\x)
    \coloneqq
    \begin{bmatrix*}[l]
        A_d(\x) & B_d^\T(\x)
        \\[1ex]
        B_d(\x) & A_d(\x)
    \end{bmatrix*}
    =
    \begin{bmatrix*}[r]
        A_d(\x) & -B_d(\x)
        \\[1ex]
        B_d(\x) &  A_d(\x)
    \end{bmatrix*}
    \MXkern.
\end{equation}
Let
\begin{equation*}
    U
    \coloneqq
    \frac{\sqrt{2}}{2}
    \begin{bmatrix}
        I_d     & \ii I_d
        \\
        \ii I_d & I_d
    \end{bmatrix}
    \MXkern,
\end{equation*}
where $I_d$ is the identity matrix of size~$d$.
After the change of coordinates
\begin{equation*}
    \overline{U} A_{2d} U^\T
    =
    \begin{bmatrix*}
        A + \ii B & 0
        \\
        0         & A - \ii B
    \end{bmatrix*}
    =
    \begin{bmatrix*}
        M_d & 0
        \\
        0   & \overline{M}_d
    \end{bmatrix*}
    \MXkern,
\end{equation*}
we see that
\begin{equation*}
    \det(A_{2d}(\x))
    =
    \det\big(M_d(\x)\big) \det\big(\overline{M}_d(\x)\big)
    =
    f^2.
\end{equation*}
The polynomials $f$ and $f^2$ have the same hyperbolicity cone
with respect to $\e$,
which is the spectrahedron defined by $A_{2d}(\x)$.

\bigskip\noindent
The \emph{rank} and \emph{corank} of a point $\x \in \C\P^n$
are defined as $\rank M_d(\x)$ and $\corank M_d(\x)$, respectively.
The \locus{k} of $M_d(\x)$ is the set of points in $\C\P^n$
with rank less than or equal to $k$.
The hypersurface~$\V(f)$ is equal to the \locus{(d - 1)}.
The points $\x \in \V(f)$ with $\corank(\x) \ge 2$ with respect to $M_d(\x)$
are always singular on $V(f)$,
and they are called \emph{essential singularities}.
A point $\x$ with $\corank(\x) = 1$ is generally not singular,
but if $\x \in \Sing \V(f)$,
then $\x$ is called an \emph{accidental singularity}.
The multiplicity of a point $\x \in \V(f)$ 
is greater than or equal to its corank.
Since the \locus{(d - 2)} is given
by the vanishing of the $(d - 1) \times (d - 1)$-minors of $M_d(\x)$,
the singular locus of $\V(f)$ is at least $(n - 4)$-dimensional,
and it is at least $(n - 3)$-dimensional if $\V(f)$ is a symmetroid.
Moreover,
if $\V(f)$ is a generic symmetroid of degree~$d$,
then $\Sing \V(f)$ has degree~$\binom{d + 1}{3}$
and contains no accidental singularities \cites[Proposition~1.5]{Pio06}[420]{Sal65}.

We restrict the attention to quartic determinantal surfaces in $\P^3$\Pnkern.
A generic determinantal surface is smooth,
while a generic quartic symmetroid has ten \point{2}s that are nodes,
that is,
isolated quadratic singularities.
A nodal quartic symmetroid is called \emph{transversal}
if it has ten \node{2}s and no further singularities.
The study of quartic symmetroids originated with Cayley \cite{Cay69}.
Recently,
real quartic symmetroids with a nonempty spectrahedron
have gained attention.
Using the global Torelli theorem for K3-surfaces,
Degtyarev and Itenberg proved the following:

\begin{theorem}[{\cite[Theorem~1.1]{DI11}}]
    \label{thm:DI}
    There exists a real transversal quartic symmetroid
    with a nonempty spectrahedron,
    having $\rho$ real nodes,
    of which $\sigma$ nodes lie on the boundary of the spectrahedron,
    if and only if $0 \le \sigma \le \rho$, both even,
    and $2 \le \rho \le 10$.
\end{theorem}

\noindent
Ottem et al.\ presented an algorithmic proof of \cref{thm:DI},
and for each pair $(\rho, \sigma)$ satisfying the inequalities in the theorem,
they gave an example of a symmetroid
with the corresponding configuration of nodes \cite{Ott+14}.
An analogue to \cref{thm:DI} for rational quartic symmetroids
is proven in \cite[Theorem~1.7]{HR18}.

We prove the following bound on the number of isolated \point{2}s
of a determinantal quartic surface in \cref{sec:bound}:

\begin{theorem}
    \label{thm:essential-sing}
    Let $S_4 \subset \P^3$ be a determinantal quartic surface
    with only isolated, simple singularities.
    Let $\eta$ be the number of essential singularities of $S_4$.
    Then $\eta \le 8$,
    unless $S_4$ is a symmetroid,
    in which case $\eta = 10$.
\end{theorem}

\noindent
We expect that \cref{thm:essential-sing} is well-known,
but we have not been able to find a reference in the literature.
As remarked above,
a definite Hermitian determinantal representation
gives rise to a spectrahedron.
The real singularities may lie on or off the spectrahedron.
We surmise a generalisation of \cref{thm:DI}
to Hermitian representations:

\begin{conjecture}
    \label{conj:hermitian}
    Suppose that $\V(f) \subset \R\P^3$ is a real quartic surface,
    where $f$ admits
    a definite Hermitian determinantal representation~$M_4(\x)$.
    Assume that the complex surface~$V_\C(f) \subset \C\P^3$
    has $\eta$ isolated nodes,
    all of which are essential nodes with respect to $M_4(\x)$,
    and that $\V(f)$ has $\rho$ real nodes,
    of which
    $\sigma$ real nodes lie on the spectrahedron defined by $M_4(\x)$.
    Then $f$ exists if and only if
    $0 \le \sigma \le \rho \le \eta \le 8$
    and
    $\rho \equiv \eta \pmod{2}$,
    or $0 \le \sigma \le \rho \le \eta = 10$,
    $\rho \ge 2$ and $\sigma \equiv \rho \equiv \eta \pmod{2}$.
\end{conjecture}

\noindent
In \cref{sec:real-sing},
we consider the real singularities of quartic surfaces with a Hermitian determinantal representation.
In particular,
\cref{sec:definite-rep}
explains why we expect fewer restrictions on
$\rho$ and $\sigma$ for $0 \le \eta \le 8$ than for $\eta = 10$ in \cref{conj:hermitian}.
After that,
we describe our strategy for finding examples of surfaces
with a given triple~$(\eta, \rho, \sigma)$.
\cref{tab:missing-cases} shows the progress towards proving \cref{conj:hermitian}.
The existence of all cases for $\eta = 10$ is given by \cref{thm:DI},
and explicit examples are given in \cite{Ott+14}.
The known examples for $0 \le \eta \le 8$ are listed in \cref{sec:hermitian-examples}.

\begin{longtable}{@{}
                  >{\(}m{0.05\textwidth}<{\)}
                  m{0.91\textwidth}
                  @{}}

    \caption{Progress on \cref{conj:hermitian}.}
    \label{tab:missing-cases}
    \\*[-1ex]
    \toprule
    \boldsymbol{\eta}
    &
    \(\boldsymbol{(\rho, \sigma)}\)
    \\*
    \midrule
    \endfirsthead

    \toprule
    \boldsymbol{\eta}
    &
    \(\boldsymbol{(\rho, \sigma)}\)
    \\*
    \midrule
    \endhead

    \addlinespace
    \multicolumn{2}{@{}l}{Continued on the next page.}
    \endfoot
 
    \bottomrule
    \endlastfoot

        0
        &
        Known examples:
        \hspace{0.1em}
        $(0, 0)$
        \vskip1ex
        Missing examples:
        \hspace{0.01em}
        None
        \\*
        \midrule

        1
        &
        Known examples:
        \hspace{0.1em}
        $(1, 1)$,
        $(1, 0)$
        \vskip1ex
        Missing examples:
        \hspace{0.01em}
        None
        \\*
        \midrule

        2
        &
        Known examples:
        \hspace{0.1em}
        $(2, 2)$,
        $(2, 1)$,
        $(2, 0)$,
        $(0, 0)$
        \vskip1ex
        Missing examples:
        \hspace{0.01em}
        None
        \\*
        \midrule

        3
        &
        Known examples:
        \hspace{0.1em}
        $(3, 3)$,
        $(3, 2)$,
        $(3, 1)$,
        $(3, 0)$,
        $(1, 1)$,
        $(1, 0)$
        \vskip1ex
        Missing examples:
        \hspace{0.01em}
        None
        \\*
        \midrule

        4
        &
        Known examples:
        \hspace{0.1em}
        $(4, 4)$,
        $(4, 3)$,
        $(4, 2)$,
        $(4, 1)$,
        $(4, 0)$,
        $(2, 2)$,
        $(2, 1)$,
        $(0, 0)$
        \vskip1ex
        Missing examples:
        \hspace{-0.15em}
        $(2, 0)$
        \\*
        \midrule

        5
        &
        Known examples:
        \hspace{0.1em}
        $(5, 3)$,
        $(5, 2)$,
        $(5, 1)$
        \vskip1ex
        Missing examples:
        $(5, 5)$,
        $(5, 4)$,
        $(5, 0)$,
        $(3, 3)$,
        $(3, 2)$,
        $(3, 1)$,
        $(3, 0)$,
        $(1, 1)$,
        $(1, 0)$
        \\*
        \midrule

        6
        &
        Known examples:
        \hspace{0.1em}
        $(6, 4)$,
        $(6, 3)$
        \vskip1ex
        Missing examples:
        $(6, 6)$,
        $(6, 5)$,
        $(6, 2)$,
        $(6, 1)$,
        $(6, 0)$,
        $(4, 4)$,
        $(4, 3)$,
        $(4, 2)$,
        $(4, 1)$,
        $(4, 0)$,
        $(2, 2)$,
        $(2, 1)$,
        $(2, 0)$,
        $(0, 0)$
        \\*
        \midrule

        7
        &
        Known examples:
        \hspace{0.1em}
        $(7, 5)$,
        $(7, 4)$
        \vskip1ex
        Missing examples:
        $(7, 7)$,
        $(7, 6)$,
        $(7, 3)$,
        $(7, 2)$,
        $(7, 1)$,
        $(7, 0)$,
        $(5, 5)$,
        $(5, 4)$,
        $(5, 3)$,
        $(5, 2)$,
        $(5, 1)$,
        $(5, 0)$,
        $(3, 3)$,
        $(3, 2)$,
        $(3, 1)$,
        $(3, 0)$,
        $(1, 1)$,
        $(1, 0)$
        \\*
        \midrule

        8
        &
        Known examples:
        \hspace{0.1em}
        $(8, 5)$,
        $(8, 4)$,
        $(6, 4)$
        \vskip1ex
        Missing examples:
        $(8, 8)$,
        $(8, 7)$,
        $(8, 6)$,
        $(8, 3)$,
        $(8, 2)$,
        $(8, 1)$,
        $(8, 0)$,
        $(6, 6)$,
        $(6, 5)$,
        $(6, 3)$,
        $(6, 2)$,
        $(6, 1)$,
        $(6, 0)$,
        $(4, 4)$,
        $(4, 3)$,
        $(4, 2)$,
        $(4, 1)$,
        $(4, 0)$,
        $(2, 2)$,
        $(2, 1)$,
        $(2, 0)$,
        $(0, 0)$
        \\*
        \midrule

        10
        &
        Known examples:
        \hspace{0.1em}
        $(10, 10)$,
        $(10, 8)$,
        $(10, 6)$,
        $(10, 4)$,
        $(10, 2)$,
        $(10, 0)$,
        $(8, 8)$,
        $(8, 6)$,
        $(8, 4)$,
        $(8, 2)$,
        $(8, 0)$,
        $(6, 6)$,
        $(6, 4)$,
        $(6, 2)$,
        $(6, 0)$,
        $(4, 4)$,
        $(4, 2)$,
        $(4, 0)$,
        $(2, 2)$,
        $(2, 0)$
        \vskip1ex
        Missing examples:
        \hspace{0.01em}
        None
\end{longtable}

\section{Essential singularities on determinantal quartic surfaces}
\label{sec:bound}

A quartic surface~$\V(f) \subset \P^3$ with only isolated singularities,
can have zero to sixteen nodes.
If $f$ has a symmetric determinantal representation,
exactly ten of the nodes
\dash counted with multiplicity \dash
are essential singularities.
For each of the sixteen nodes on a Kummer surface~$\V(f)$,
there exists a symmetric determinantal representation of $f$
such that the node is essential in that representation \cite[597]{Ott+14}.
On the other hand,
\cite[Article~9]{Jes16} describes a quartic symmetroid with eleven nodes,
where one of the nodes is an accidental singularity
in every symmetric determinantal representation.
It is natural to ask how many nodes can be essential singularities 
when we consider \emph{nonsymmetric} determinantal representations.
We show that the maximum number of isolated, essential singularities
is obtained precisely with a symmetric determinantal representation.

Determinantal quartic surfaces are characterised
by containing a projectively normal sextic curve with genus~$3$ \cite{Sch81}.
Coble noted that this follows because the Picard group of a general quartic surface
is generated by a plane section \cite[39]{Cob82}.
Given a determinantal representation
\begin{equation*}
    M_4(\x)
    \coloneqq
    \begin{bmatrix}
        m_{00} & m_{10} & m_{20} & m_{30}
        \\
        m_{01} & m_{11} & m_{21} & m_{31}
        \\
        m_{02} & m_{12} & m_{22} & m_{32}
        \\
        m_{03} & m_{13} & m_{23} & m_{33}
    \end{bmatrix}
\end{equation*}
of a quartic surface $S_4 \coloneqq \V(\det (M_4(\x)))$,
two families of genus~$3$ sextics on $S_4$ can be described as follows:
The $(3 \times 3)$-minors of a $(4 \times 3)$-submatrix of $M_4(\x)$
define a sextic curve of genus~$3$;
the four curves obtained by the different $(4 \times 3)$-submatrices span the first family~$\mathcal{C}_1$.
The second family~$\mathcal{C}_2$ is defined similarly
by replacing the $(4 \times 3)$-submatrices with $(3 \times 4)$-submatrices.
If $M_4(\x)$ is symmetric,
$\mathcal{C}_1$ and $\mathcal{C}_2$ coincide;
if the representation is Hermitian,
the two families are complex conjugates.
Because these curves are given by the vanishing of some $(3 \times 3)$-minors,
they contain the \locus{2} of $S_4$.
That is,
they contain the set of essential singularities on $S_4$.
Moreover,
for each curve $C_1 \in \mathcal{C}_1$ there is a curve $C_2 \in \mathcal{C}_2$ 
such that $C_1 \cup C_2$ is the complete intersection of $S_4$ and a cubic surface $S_3$.
In particular,
the union of the curves defined by the $(3 \times 3)$-minors of
\begin{equation*}
    A_1
    \coloneqq
    \begin{bmatrix}
        m_{00} & m_{10} & m_{20} 
        \\
        m_{01} & m_{11} & m_{21} 
        \\
        m_{02} & m_{12} & m_{22} 
        \\
        m_{03} & m_{13} & m_{23}
    \end{bmatrix}
    \qquad
    \text{and}
    \qquad
    A_2
    \coloneqq
    \begin{bmatrix}
        m_{00} & m_{10} & m_{20} & m_{30}
        \\
        m_{01} & m_{11} & m_{21} & m_{31}
        \\
        m_{02} & m_{12} & m_{22} & m_{32}
    \end{bmatrix}
\end{equation*}
is the intersection of $S_4$ and the surface~$S_3$ defined by their common $(3 \times 3)$-minor
\begin{equation*}
    \begin{vmatrix}
        m_{00} & m_{10} & m_{20} 
        \\
        m_{01} & m_{11} & m_{21} 
        \\
        m_{02} & m_{12} & m_{22} 
    \end{vmatrix}
    =
    0.
\end{equation*}

The curves in $\mathcal{C}_1$ and $\mathcal{C}_2$ are nonhyperelliptic.
Indeed,
a sextic curve of genus~$3$ in $\P^3$ is nonhyperelliptic
if and only if it is projectively normal \cite[Exercise~4.10]{Dol12}.
We give an elementary argument showing that the curve~$C$
defined by $A_1$, is nonhyperelliptic:
Suppose that the entries $m_{ij}$ in $A_1$
are linear forms in the variables $x_0$, $x_1$, $x_2$ and $x_3$.
Note that $C$ is the solution set to the equation
\begin{equation}
    \label{eq:matrix-equality}
    A_1
    \begin{bmatrix}
        y_0 \\ y_1 \\ y_2
    \end{bmatrix}
    =
    \begin{bmatrix}
        0 \\ 0 \\ 0 \\ 0
    \end{bmatrix}
    \MXkern.
\end{equation}
We can rewrite \eqref{eq:matrix-equality} as
\begin{equation}
    \label{eq:matrix-equality-rewritten}
    A'_1
    \begin{bmatrix}
        x_0 \\ x_1 \\ x_2 \\ x_3
    \end{bmatrix}
    =
    \begin{bmatrix}
        0 \\ 0 \\ 0 \\ 0
    \end{bmatrix}
    \MXkern,
\end{equation}
where $A'_1$ is a $(4 \times 4)$-matrix
with linear entries in $y_0$, $y_1$ and $y_3$.
Both \eqref{eq:matrix-equality} and \eqref{eq:matrix-equality-rewritten}
define the same curve~$K$ in $\P^2 \times \P^3$\Pnkern.
Then $C$ is the projection of $K$ to $\P^3$\Pnkern,
and the curve~$C'$ given by $\det(A'_1) = 0$ is the projection to $\P^2$\Pnkern.
The curve~$C'$ is a smooth planar quartic curve,
hence nonhyperelliptic \cite[Example~IV.5.2.1]{Har77}.
It follows that $C$ is nonhyperelliptic as well.

We are now ready to prove the bound on the possible number
of essential singularities on a quartic surface.

\begin{proof}[Proof of \cref{thm:essential-sing}.]
    Without loss of generality,
    we may assume that $S_4$ has only essential singularities,
    $P_1, \ldots, P_\eta$.
    Since $S_4$ is determinantal,
    there are smooth, sextic curves $C_1, C_2 \subset S_4$ of genus~$3$
    passing through $P_1, \ldots, P_\eta$.
    Moreover,
    we may assume that $C_1 \cup C_2 = S_4 \cap S_3$ for some cubic surface~$S_3$.

    Let $\pi \colon \widetilde{S}_4 \to S_4$ be the blow-up of $S_4$ at $P_1, \ldots, P_\eta$.
    Then $\widetilde{S}_4$ is a smooth K3-surface.
    The exceptional divisor~$E_i$ over $P_i$
    is a $(-2)$-curve satisfying $h \cdot E_i = 0$,
    where $h$ is the class of the preimage of a plane section of $S_4$.
    Because $C_1$ passes through $P_i$,
    the strict transform~$\widetilde{C}_1$ satisfies $\widetilde{C}_1 \cdot E_i = 1$.
    Furthermore,
    the adjunction formula gives
    \begin{equation*}
        \label{eq:adjunction}
        \widetilde{C}_1 \cdot \widetilde{C}_1
        =
        2g_{\widetilde{C}_1} - 2
        =
        2 \cdot 3 - 2
        =
        4,
    \end{equation*}
    since the canonical divisor on $\widetilde{S}_4$ is trivial.
    It follows from $C_1 \cup C_2 = S_4 \cap S_3$
    that the total transform~$\pi^{-1} (C_1 \cup C_2) = 3h$.
    The curve $C_1 \cup C_2$ is double at $P_1, \ldots, P_\eta$,
    so
    $
        \pi^{-1} (C_1 \cup C_2)
        =
        \widetilde{C}_1 + \widetilde{C}_2 - \sum_{i = 1}^\eta E_i
    $.
    Thus
    \begin{equation}
        \label{eq:sum-of-strict-transforms}
        \widetilde{C}_1 + \widetilde{C}_2
        =
        \pi^{-1} (C_1 \cup C_2) - \sum_{i = 1}^\eta E_i
        =
        3h - \sum_{i = 1}^\eta E_i.
    \end{equation}
    We intersect both sides of \eqref{eq:sum-of-strict-transforms} 
    with $\widetilde{C}_1$:
    \begin{equation}
        \label{eq:transforms-intersected}
        \widetilde{C}_1\cdot \Big(\widetilde{C}_1 + \widetilde{C}_2 \Big)
        =
        \widetilde{C}_1 \cdot \Bigg(3h - \sum_{i = 1}^\eta E_i\Bigg).
    \end{equation}
    By using the facts above and that
    $\widetilde{C}_1 \cdot h = \deg\big(\widetilde{C}_1\big) = 6$,
    \cref{eq:transforms-intersected} yields
    $\widetilde{C}_1 \cdot \widetilde{C}_2 = 14 - \eta$.

    The linear system $\abs[\big]{\widetilde{C}_1}$
    gives rise to a morphism $\varphi \colon \widetilde{S}_4 \to \P^3$.
    Because $\widetilde{C}_1$ is a nonhyperelliptic curve of genus~$3$ with $\widetilde{C}_1 \cdot \widetilde{C}_1 = 4$,
    it is mapped to a plane curve~$C'_1 \coloneqq \varphi\big(\widetilde{C}_1\big)$ of degree~$4$.
    Thus the image $S'_4 \coloneqq \varphi \big( \widetilde{S}_4 \big)$ is a quartic surface.
    Assume first that $S_4$ is not a symmetroid,
    so $C_1 \neq C_2$.
    Then $C'_2 \coloneqq \varphi\big(\widetilde{C}_2\big)$ is not a plane section.
    Because $C'_2$ is nonhyperelliptic and spans $\P^3$\Pnkern,
    we have $\deg(C'_2) \ge 6$.
    Since $C'_1$ is a plane section,
    we have $\deg(C'_2) = \widetilde{C}_1 \cdot \widetilde{C}_2$.
    Therefore
    \begin{equation}
        \label{eq:intersection-products}
        14 - \eta
        =
        \widetilde{C}_1 \cdot \widetilde{C}_2
        \ge
        6,
    \end{equation}
    so we get $\eta \le 8$.

    Assume now that $S_4$ is a symmetroid,
    so $C_1 = C_2$.
    Then $\widetilde{C}_1 \cdot \widetilde{C}_2 = \widetilde{C}_1 \cdot \widetilde{C}_1 = 4$.
    Adjusting for this in the right-hand side of \eqref{eq:intersection-products} gives $\eta = 10$. 
    Hence we have recovered the well-known fact
    that a quartic symmetroid with only isolated singularities has ten \point{2}s. 
\end{proof}

\section{Real singularities on quartics
         with a Hermitian determinantal representation}
\label{sec:real-sing}

Consider $\C\P^{15}$
as the projectivisation of the vector space over $\C$
spanned by Hermitian $(4 \times 4)$-matrices.
Then the \locus{2}~$X_2$ of $\C\P^{15}$ is given by the vanishing
of the $(3 \times 3)$-minors of a general Hermitian $(4 \times 4)$-matrix;
it is an elevenfold of degree~$20$.
A quartic surface $S_4 \subset \P^3$ with a Hermitian determinantal representation
corresponds to a linear $3$-space~$H \subset \C\P^{15}$\Pnkern,
and essential singularities on $S_4$ corresponds to the intersection
of $H$ with $X_2$.
This is helpful for finding examples
of Hermitian determinantal representations with a specific singular locus.

In \cite{HR18, Hel19, Hel20},
symmetroids are studied via the associated quadratic form.
We wish to use this technique to determine the real part of $X_2$.
Let $\y \coloneqq [y_0, y_1, y_2, y_3]$.
The form~$h_{M_4}(\y) \coloneqq \overline{\y} \mkern-1mu M_4 \y^\T$
associated to a Hermitian $(4 \times 4)$-matrix~$M_4$
is not polynomial in $y_0$, $y_1$, $y_2$ and $y_3$.
However,
we can associate a quadratic form
to the symmetric matrix~$A_8$ from \eqref{eq:associated-sym}.
Incidentally,
$h_{M_4} = \y'\mkern-2mu A_8 \y'^\T$,
where
\begin{equation*}
    \y'
    \coloneqq
    [\Re(y_0), \Re(y_1), \Re(y_2), \Re(y_3),
     \Im(y_0), \Im(y_1), \Im(y_2), \Im(y_3)].
\end{equation*}
Since $M_4$ is Hermitian,
it is unitarily diagonalisable.
If $M_4$ is diagonal,
then so is $A_8$.
From this we see that 
$M_4$ and $A_8$ have the same eigenvalues~$\lambda_i$,
but the algebraic multiplicity $\mu_{A_8}(\lambda_i)$ is $2\mu_{M_4}(\lambda_i)$.
Hence $\rank A_8 = 2 \rank(M_4)$.
In addition,
if $M_4(\x)$ is a Hermitian determinantal representation
of a real polynomial~$f$,
then we can define its spectrahedron in terms of $M_4(\x)$ only,
because
\begin{equation*}
    \big\{
        \x \in \R\P^3
        \mid
        M_4(\x) \text{ is semidefinite}
    \big\}
    =
    \big\{
        \x \in \R\P^3
        \mid
        A_8(\x) \text{ is semidefinite}
    \big\}.
\end{equation*}

\newcommand{\kerncolon}{\mkern-1mu:\mkern-1mu}
The construction of $A_8$ from $M_4$ allows us to view the point
\begin{equation*}
    \x
    \coloneqq
    [x_{00} \kerncolon x_{01} \kerncolon x_{02} \kerncolon x_{03}
            \kerncolon x_{11} \kerncolon x_{12} \kerncolon x_{13}
                              \kerncolon x_{22} \kerncolon x_{33}
            \kerncolon y_{01} \kerncolon y_{02} \kerncolon y_{03}
                              \kerncolon y_{12} \kerncolon y_{13}
                              \kerncolon y_{23}]
\end{equation*}
in $\C\P^{15}$ as the symmetric $(8 \times 8)$-matrix
\begin{equation*}
    \label{eq:8x8}
    A_8(\x)
    \coloneqq
    \left[
    \begin{array}{@{}
                 c@{\hspace{1ex}}
                 c@{\hspace{1ex}}
                 c@{\hspace{1ex}}
                 c@{\hspace{3ex}}
                 c@{\hspace{1ex}}
                 c@{\hspace{1ex}}
                 c@{\hspace{1ex}}
                 c@{}}
         x_{00} &  x_{01} &  x_{02} &  x_{03} &      
              0 & -y_{01} & -y_{02} & -y_{03}
        \\
         x_{01} &  x_{11} &  x_{12} &  x_{13} &
         y_{01} &       0 & -y_{12} & -y_{13} 
        \\
         x_{02} &  x_{12} &  x_{22} &  x_{23} &
         y_{02} &  y_{12} &       0 & -y_{23}
        \\
         x_{03} &  x_{13} &  x_{23} &  x_{33} &
         y_{03} &  y_{13} &  y_{23} &       0
        \\[2ex]
              0 &  y_{01} &  y_{02} &  y_{03} &
         x_{00} &  x_{01} &  x_{02} &  x_{03}
        \\
        -y_{01} &       0 &  y_{12} &  y_{13} &
         x_{01} &  x_{11} &  x_{12} &  x_{13}
        \\
        -y_{02} & -y_{12} &       0 &  y_{23} &
         x_{02} &  x_{12} &  x_{22} &  x_{23}
        \\
        -y_{03} & -y_{13} & -y_{23} &      0  &
         x_{03} &  x_{13} &  x_{23} &  x_{33}
    \end{array}
    \right]
    \MXkern.
\end{equation*}
Let $\y \coloneqq [y_0, \ldots, y_7]$.
Denote by $Q_{A_8}(\x)$ the quadric $\V\big(\y \mkern-1mu A_8(\x) \y^\T\big) \subset \C\P^7$ 
associated to $A_8(\x)$,
and let
\[
    W
    \coloneqq
    \Big\{
        Q_{A_8(\x)} \mid \x \in \C\P^{15}
    \Big\}.
\]
After a \cite{Macaulay2} calculation,
we find that the baselocus~$\Bl(W) \subset \C\P^7$ of $W$
consists of the two disjoint, complex conjugate $3$-spaces
$H_{\Bl} \coloneqq \V(l_0, l_1, l_2, l_3)$ and
\(
    \overline{H}_{\Bl}
    \coloneqq
    \V\big(
        \overline{l}_0, \overline{l}_1,\overline{l}_2, \overline{l}_3
      \big)
\),
where $l_j$ is the linear form~$y_j + \ii y_{j + 4}$ for $j = 0, 1, 2, 3$.
We summarise the correspondence between $M_4$ and $A_8$:

\begin{lemma}
    \label{lem:hermitian-quadrics}
    There is an isomorphism between the
    projectivisation of the vector space over $\C$
spanned by Hermitian $(4 \times 4)$-matrices of rank~$k$,
    and
    the space of \quadric{2k}s in $\C\P^7$
    passing through two given disjoint $3$-spaces.
\end{lemma}

\noindent
A point in $X_2$ corresponds to a \quadric{4} in $W$\Wkern.
We will now count the \quadric{4}s in $W$\Wkern.
Let $Q$ be a \quadric{4} in $W$\Wkern.
By the rank-nullity theorem,
$\Sing(Q)$ is a $3$-space.
Since $Q$ contains $\Bl(W)$ and is irreducible,
it follows that $\Sing(Q)$ intersects $H_{\Bl}$ and $\overline{H}_{\Bl}$
in a line each.
A $3$-space~$H$ which intersects $H_{\Bl}$ and $\overline{H}_{\Bl}$
in a line each,
is spanned by those lines,
since $H_{\Bl}$ and $\overline{H}_{\Bl}$ are disjoint.
The Grassmannian~$\G(1, 3)$ of lines in a $3$-space is $4$-dimensional.
Hence there is an $8$-dimensional space of $3$-spaces
that are spanned by a line in $H_{\Bl}$ and a line in $\overline{H}_{\Bl}$.
The $3$-space $H$ is the singular locus
of a web $W_H \subset W$ of quadrics.
This can be seen by projecting $\P^7 \to \P^3$
with $H$ as projection centre.
The $3$-spaces $H_{\Bl}$ and $\overline{H}_{\Bl}$ are projected
onto two skew lines $L_{\Bl}, \overline{L}_{\Bl} \subset \P^3$\Pnkern.
A \quadric{4}~$Q \subset \P^7$ singular at $H$
and containing $H_{\Bl}$ and $\overline{H}_{\Bl}$
is projected onto a quadric in $\P^3$
containing $L_{\Bl}$ and $\overline{L}_{\Bl}$.
There is a web of quadrics in $\P^3$ containing $L_{\Bl}$ and $\overline{L}_{\Bl}$.
In total,
we get that the space of \quadric{4}s in $W$ has dimension $8 + 3 = 11$.
Since there is a bijection between this space and $X_2$,
this is as expected.

A real point in $X_2$ corresponds to a real \quadric{4} $Q$ in $W$\Wkern.
Then $H \coloneqq \Sing(Q)$ is real.
Thus $H$ intersects $H_{\Bl}$ and $\overline{H}_{\Bl}$
in two complex conjugate lines, $L$, $\overline{L}$.
Hence $H_L \coloneqq \linspan[\big]{L, \overline{L}}$
is the unique real $3$-space which contains $L$
and is the singular locus of a \quadric{4} in $W$\Wkern.
As noted above,
there is a web of quadrics singular at $H_L$.
Since the Grassmannian~$\G(1, 3)$ of lines in $H_{\Bl}$
is $4$-dimensional,
we get that the locus of real \quadric{4}s in $W$ has dimension $3 + 4 = 7$.
Hence the real part of the \locus{2}~$X_2$ is $7$-dimensional.

\begin{proposition}
    The \locus{2}~$X_2$ of $\C\P^{15}$\Pnkern,
    the projectivisation of the vector space over $\C$
spanned by Hermitian $(4 \times 4)$-matrices,
    is $11$-dimensional.
    The real part of $X_2$ is $7$-dimensional.
\end{proposition}

\noindent
We derive a parameterisation of the real \quadric{4}s in $W$\Wkern,
which in turn corresponds to a parameterisation of the real part of $X_2$.
A line~$L$ in $H_{\Bl}$ is given by
\begin{equation*}
    \V
    \Bigg(
        y_0 + \ii y_4,
        y_1 + \ii y_5,
        y_2 + \ii y_6,
        y_3 + \ii y_7,
        \sum_{j = 0}^{7} (a_j + \ii b_j)y_j,
        \sum_{j = 0}^{7} (c_j + \ii d_j)y_j
    \Bigg)
\end{equation*}
for $a_j, b_j, c_j, d_j \in \R$.
Since $H_L$ is the only real $3$-space containing $L$,
we deduce that $H_L = V(\ell_0, \ell_1, \ell_2, \ell_3)$,
where
\begin{align*}
    \ell_0
    \coloneqq
    &\sum_{j = 0}^{3}
    \Big(
        (a_{j} - b_{j + 4})y_{j \phantom{+4}}
        +
        (a_{j + 4} + b_{j})y_{j + 4}
    \Big),
    \\
    \ell_1
    \coloneqq
    &\sum_{j = 0}^{3}
    \Big(
        (a_{j} - b_{j + 4})y_{j + 4}
        -
        (a_{j + 4} + b_{j})y_{j}
    \Big),
    \\
    \ell_2
    \coloneqq
    &\sum_{j = 0}^{3}
    \Big(
        (c_{j} - d_{j + 4})y_{j \phantom{+4}}
        +
        (c_{j + 4} + d_{j})y_{j + 4}
    \Big),
    \\
    \ell_3
    \coloneqq
    &\sum_{j = 0}^{3}
    \Big(
        (c_{j} - d_{j + 4})y_{j + 4}
        -
        (c_{j + 4} + d_{j})y_{j}
    \Big).
\end{align*}
The quadrics that are singular at $H_L$ are given by quadratic polynomials
in $\ell_0$, $\ell_1$, $\ell_2$, $\ell_3$.
Moreover,
to be contained in $W$\Wkern,
the quadrics must contain $H_{\Bl}$ and $\overline{H}_{\Bl}$.
We compute that the quadrics in $W$ that are singular at $H_L$,
generate the ideal
\begin{equation*}
    I
    \coloneqq
    \linspan[\big]
    {
        \ell_0^2  + \ell_1^2,
        \ell_0 \ell_2 + \ell_1 \ell_3,
        \ell_0 \ell_3 - \ell_1 \ell_2,
        \ell_2^2  + \ell_3^2
    }.
\end{equation*}
Hence a real \quadric{4} in $W$ is on the form
\begin{equation}
    \label{eq:parameterisation}
    \V
    \Big(
        a\big(\ell_0^2  + \ell_1^2\big)
        +
        b(\ell_0 \ell_2 + \ell_1 \ell_3)
        +
        c(\ell_0 \ell_3 - \ell_1 \ell_2)
        +
        d\big(\ell_2^2  + \ell_3^2\big)
    \Big)
\end{equation}
for $a, b, c, d \in \R$.
It would be interesting to have equations for the real part of $X_2$,
but deriving them from \eqref{eq:parameterisation}
has been too computationally demanding.

\subsection{Definite representations}
\label{sec:definite-rep}

Let $\V(f) \subset \R\P^3$ be a quartic surface
with a definite Hermitian determinantal representation~$M_4(\x)$.
Then $M_4(\x)$ gives rise to the spectrahedron
\begin{equation*}
    \big\{
        \x \in \R\P^3
        \mid
        M_4(\x) \text{ is semidefinite}
    \big\}.
\end{equation*}
If $M_4(\x)$ is semidefinite but not definite at a point $\x$,
then $\x \in \V(f)$,
and in particular,
$\x$ lies on the boundary of the spectrahedron.
Since $\V(f)$ is real,
singularities are either real or occur in complex conjugate pairs.
Hence we have $\rho \equiv \eta \pmod{2}$ in \cref{conj:hermitian}.
For $\eta < 10$,
we do not expect any further restrictions on $\rho$ and $\sigma$,
as we explain below.

Hermitian determinantal representations are subtly different
from symmetric representations,
which is why we do not have that $\sigma$ is even and $\rho \ge 2$
for all Hermitian determinantal representations with $\eta < 10$.
For instance,
the defining characteristic of a transversal symmetroid~$\V(f)$
is that the projection from one of the nodes~$P$
is ramified along the union of two cubic curves~$R_1, R_2 \subset \P^2$ \cite[200]{Cay69b}.
It is showed in \cite[Proof of Theorem~1.1, p.~600]{Ott+14}
that if $\V(f)$ has a nonempty spectrahedron~$S$
and $P \notin S$,
then the nodes on $S$ are projected onto
the intersection points of the ovals of $R_1$ and $R_2$.
These intersect in an even number of points
by the Jordan curve theorem,
hence $\sigma$ is even.
\emph{A priori},
the ramification curve for the projection from $P$
does not impose restrictions on $\sigma$
if $\V(f)$ is not a symmetroid.

In the special case of
a transversal quartic symmetroid which contains a line,
a simple reason for $\rho \neq 0$,
is that the line passes through an odd number of \point{2}s,
at least one of which must be real.
A general proof of $\rho \neq 0$ for a transversal quartic symmetroid with a nonempty spectrahedron can be found in \cite[Lemma~4.2]{Ott+14}.
Not all lines on a quartic surface with a Hermitian determinantal representation
contain an odd number of \point{2}s:

\begin{example}
    The pencil
    \begin{equation*}
        \begin{bmatrix}
            0             & x_0 - \ii x_0 & x_1 & 0
            \\
            x_0 + \ii x_0 & x_0           & x_1 & x_0 + \ii x_0
            \\
            x_1           & x_1           & x_0 & x_0
            \\
            0             & x_0 - \ii x_0 & x_0 & x_0
      \end{bmatrix}
      \end{equation*}
      has rank $3$ at all points,
      except at $[1 : 0]$ and $[0 : 1]$,
      where the rank is $2$.
\end{example}

\subsubsection{Finding examples}
In \cref{sec:hermitian-examples},
we list examples of definite Hermitian determinantal representations
with different values of $(\eta, \rho, \sigma)$.
The list is missing the sixty-four triples in \cref{tab:missing-cases}
predicted by \cref{conj:hermitian}.
Both \cite{Ott+14, HR18} present examples of symmetroids
with a nonempty spectrahedron,
that are found using random searches.
Our examples are not found in this way.
The problem with drawing random Hermitian $(4 \times 4)$-matrices $M_{4, i}$
in search of a definite determinantal representation~\eqref{eq:hermitian-representation}
with a given triple~$(\eta, \rho, \sigma)$,
is that we do not know of a matrix form for~$M_{4, i}$
that will guarantee at least $\eta \ge 4$ essential nodes.
For $\eta = 0, 1, 2, 3$
it straightforward to achieve at least $\eta$ nodes:
Let $M_{4, 0}$ be a definite matrix,
and let $\eta$ of $M_{4, 1}$, $M_{4, 2}$, $M_{4, 3}$
be \matrices{2}.
Moreover,
we can control $\sigma$ by choosing the definiteness of the \matrices{2}.
This will in general yield a representation with $(\eta, \eta, \sigma)$.
Note that $M_{4, i}$ corresponds to a real point by construction,
hence $\rho = \eta$.
To get $\rho < \eta$,
we let $M_{4, 0}$ and $M_{4, 1}$ be two \matrices{4}
that span a pencil that contains two complex conjugate \point{2}s.

For $4 \le \eta \le 8$,
our strategy has been the following:
A transversal quartic symmetroid
with a nonempty spectrahedron
corresponds to a $3$-space~$H \subset \C\P^{15}$
that intersects the \locus{2}~$X_2$ in ten points.
We try to deform $H$ whilst keeping a definite point and some nodes.
More precisely,
we choose four points $P_1$, $P_2$, $P_3$, $P_4$ in $H \cap X_2$ that span $H$\Hkern.
Let $Q \subset \P^7$ be the associated quadric at $P_1$.
There is a web of \quadric{4}s singular at $\Sing(Q)$,
which corresponds to a $3$-space~$Y_2 \subset X_2$.
We replace $P_1$ with a point $P'_1$ in $Y_2$.
The $3$-space~$H'$ spanned by $P'_1$, $P_2$, $P_3$, $P_4$
intersects $X_2$ in $\eta \ge 4$ points.
Then $H'$ corresponds to a Hermitian representation
with $\eta$ essential singularities.
There is a complete list of examples of all possible values of $\rho$ and $\sigma$ for $\eta = 10$ in \cite{Ott+14}.
We used these as a starting point
to get different values of $\rho$ and $\sigma$ for $4 \le \eta \le 8$.

\section{Examples of quartics
         with a definite Hermitian determinantal representation}
\label{sec:hermitian-examples}

Below is a list of definite Hermitian determinantal representations
that define quartic surfaces with only essential singularities.
For each matrix~$M_4(\x)$,
we specify the configuration~$(\eta, \rho, \sigma)$
of essential singularities
and a point $\e \in \R\P^3$ such that $M_4(\e)$ is definite.
We omit examples with $\eta = 10$.
A list of examples for $\eta = 10$
with all possible values for $\rho$ and $\sigma$
is found in \cite{Ott+14}.

\subsection{Zero \pointsec{2}s}

$(0, 0, 0)$:
\begin{equation*}
    \begin{bmatrix}
        x_3           & x_0           & x_2 + \ii x_1 & 0
        \\
        x_0           & x_3           & x_1           & x_1 - \ii x_2
        \\
        x_2 - \ii x_1 & x_1           & 2x_0 + x_3    & x_1 + x_2
        \\
        0             & x_1 + \ii x_2 & x_1 + x_2     & x_3
    \end{bmatrix}
\end{equation*}
$\e \coloneqq [0 : 0 : 0 : 1]$

\begin{figure}[b!htp]
    \centering
    \includegraphics[height = 0.2\textheight]{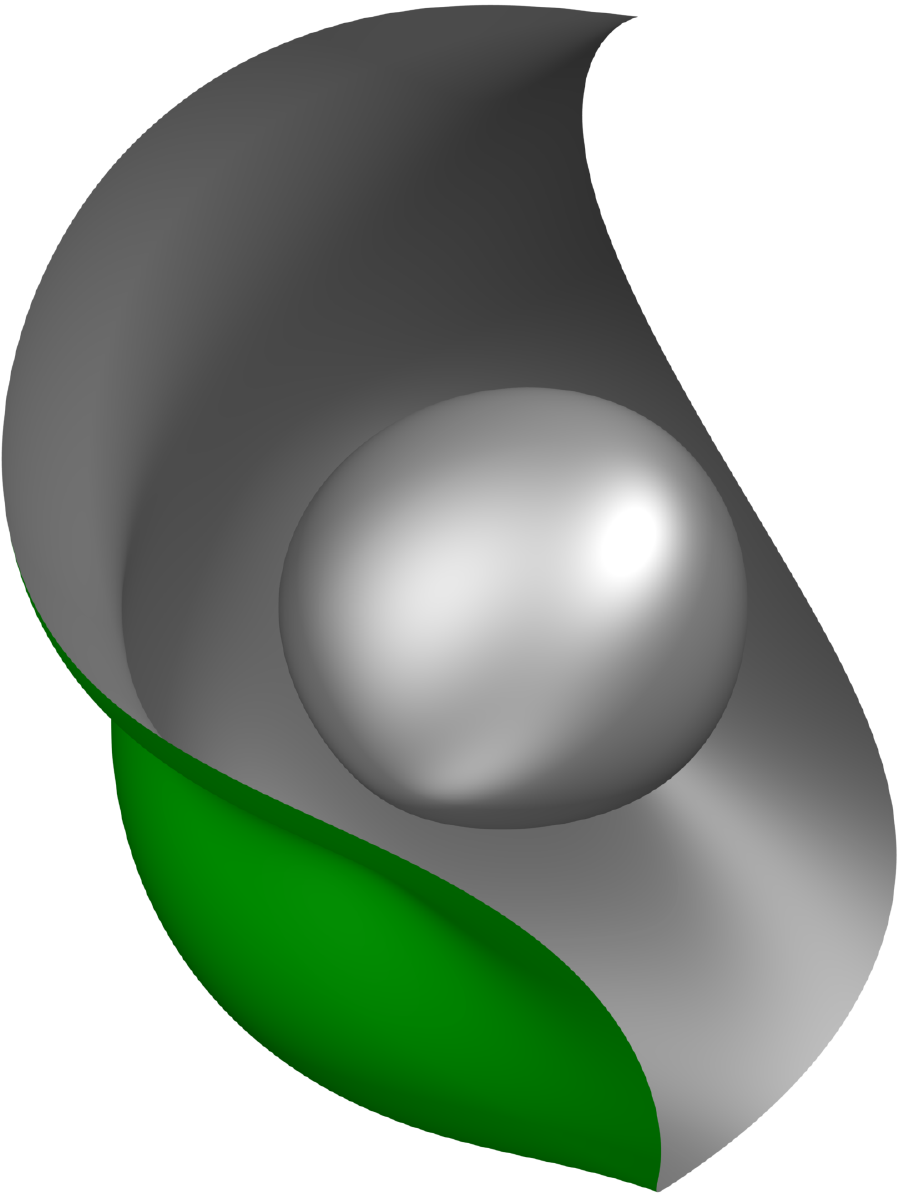}
    \addtolength{\belowcaptionskip}{-2em}
    \caption{A real, quartic surface
    with a definite Hermitian determinantal representation
    and $(\eta, \rho, \sigma) = (0, 0, 0)$.}
\end{figure}

\subsection{One \pointsec{2}}

$(1, 1, 1)$:
\begin{equation*}
    \begin{bmatrix}
        x_0 + x_3     & 0             & x_2 + \ii x_1 & 0
        \\
        0             & 2x_0 + x_3    & x_1           & x_1 - \ii x_2
        \\
        x_2 - \ii x_1 & x_1           & x_2 + x_3     & x_1 + x_2
        \\
        0             & x_1 + \ii x_2 & x_1 + x_2     & x_3
    \end{bmatrix}
\end{equation*}
$\e \coloneqq [0 : 0 : 0 : 1]$

\medskip\noindent
$(1, 1, 0)$:
\begin{equation*}
    \begin{bmatrix}
        x_3           &    x_0        & x_2 + \ii x_1 & 0
        \\
        x_0           &    x_3        &    x_1        & x_1 - \ii x_2
        \\
        x_2 - \ii x_1 &    x_1        & x_2 + x_3     & x_1 + x_2
        \\
            0         & x_1 + \ii x_2 & x_1 + x_2     & x_3
\end{bmatrix}
\end{equation*}
$\e \coloneqq [0 : 0 : 0 : 1]$

\begin{figure}[bhtp]
    \centering
    \includegraphics[height = 0.26\textheight]{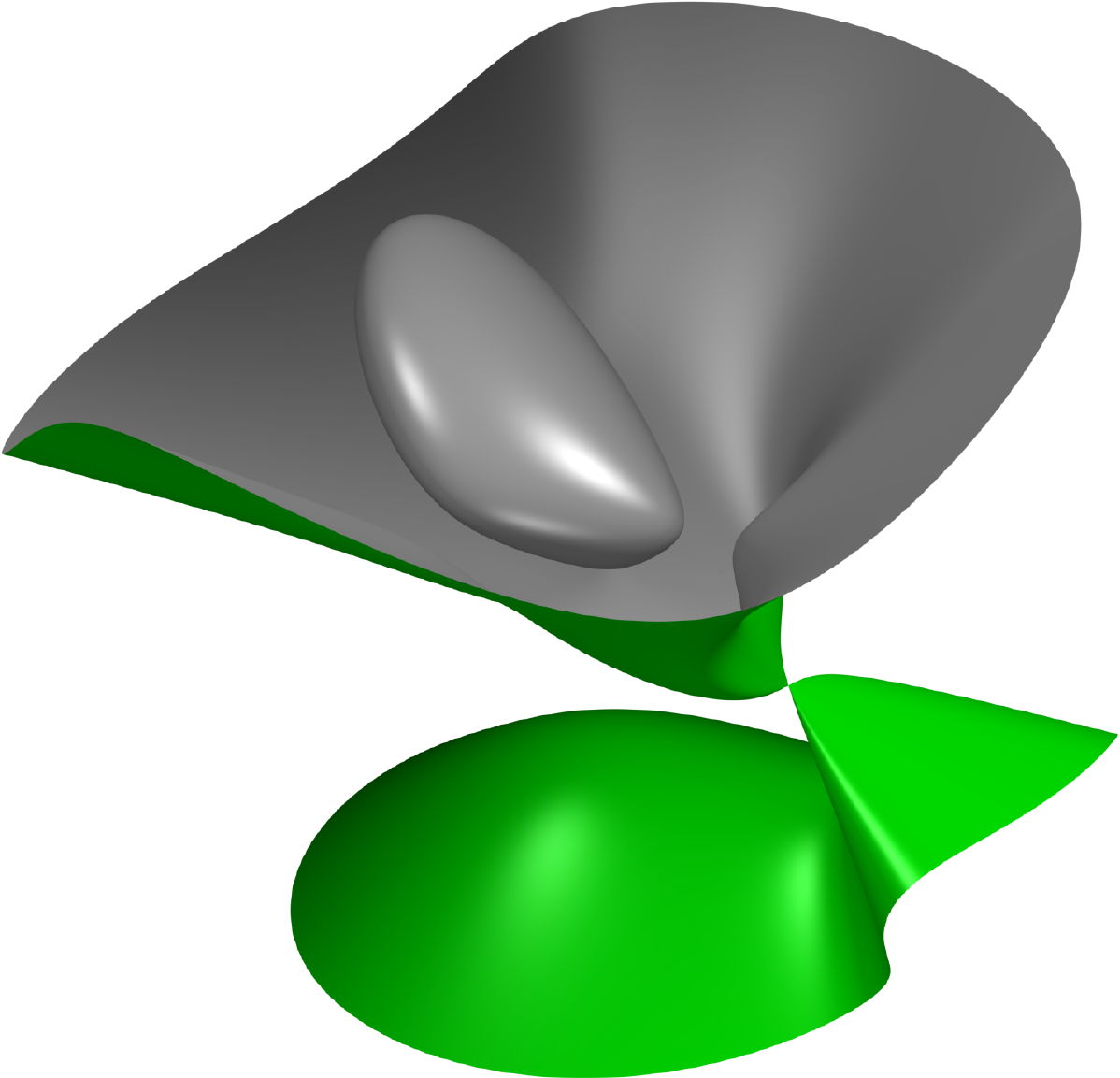}
    \caption{A real, quartic surface
    with a definite Hermitian determinantal representation
    and $(\eta, \rho, \sigma) = (1, 1, 0)$.}
\end{figure}

\subsection{Two \pointsec{2}s}

$(2, 2, 2)$:
\begin{equation*}
    \begin{bmatrix}
        x_0 + x_1 + x_3 & 0                &
        x_2 + \ii x_1   & 0
        \\
        0               & 2x_0 + x_1 + x_3 &
        x_1             & x_1 - \ii x_2
        \\
        x_2 - \ii x_1   & x_1              &
        2x_1 + x_3      & x_1 + x_2
        \\
        0               & x_1 + \ii x_2    &
        x_1 + x_2       & x_1 + x_3
    \end{bmatrix}
\end{equation*}
$\e \coloneqq [0 : 0 : 0 : 1]$

\medskip\noindent
$(2, 2, 1)$:
\begin{equation*}
    \begin{bmatrix}
        x_1 + x_3     & x_0           & x_2 + \ii x_1 & 0
        \\
        x_0           & x_1 + x_3     & x_1           & x_1 - \ii x_2
        \\
        x_2 - \ii x_1 & x_1           & 2x_1 + x_3    & x_1 + x_2
        \\
        0             & x_1 + \ii x_2 & x_1 + x_2     & x_1 + x_3
    \end{bmatrix}
\end{equation*}
$\e \coloneqq [0 : 0 : 0 : 1]$

\medskip\noindent
$(2, 2, 0)$:
\begin{equation*}
    \begin{bmatrix}
        x_3 & x_0     & x_2           & x_1
        \\
        x_0 & x_3     & 0             & -\ii x_2
        \\
        x_2 & 0       & x_3           & x_2 + \ii x_1
        \\
        x_1 & \ii x_2 & x_2 - \ii x_1 & x_3
    \end{bmatrix}
\end{equation*}
$\e \coloneqq [0 : 0 : 0 : 1]$

\medskip\noindent
$(2, 0, 0)$:
\begin{equation*}
    \begin{bmatrix}
        2x_3           & x_2           & x_2 + x_3     & x_0 + \ii x_1
        \\
        x_2            & x_3           & x_0 + \ii x_1 & 0
        \\
        x_2 + x_3      & x_0 - \ii x_1 & x_2 + x_3     & x_2
        \\
         x_0 - \ii x_1 & 0             & x_2           & x_2 + x_3
    \end{bmatrix}
\end{equation*}
$\e \coloneqq [0 : 0 : 0 : 1]$

\begin{figure}[hbtp]
    \centering
    \includegraphics[height = 0.29\textheight]{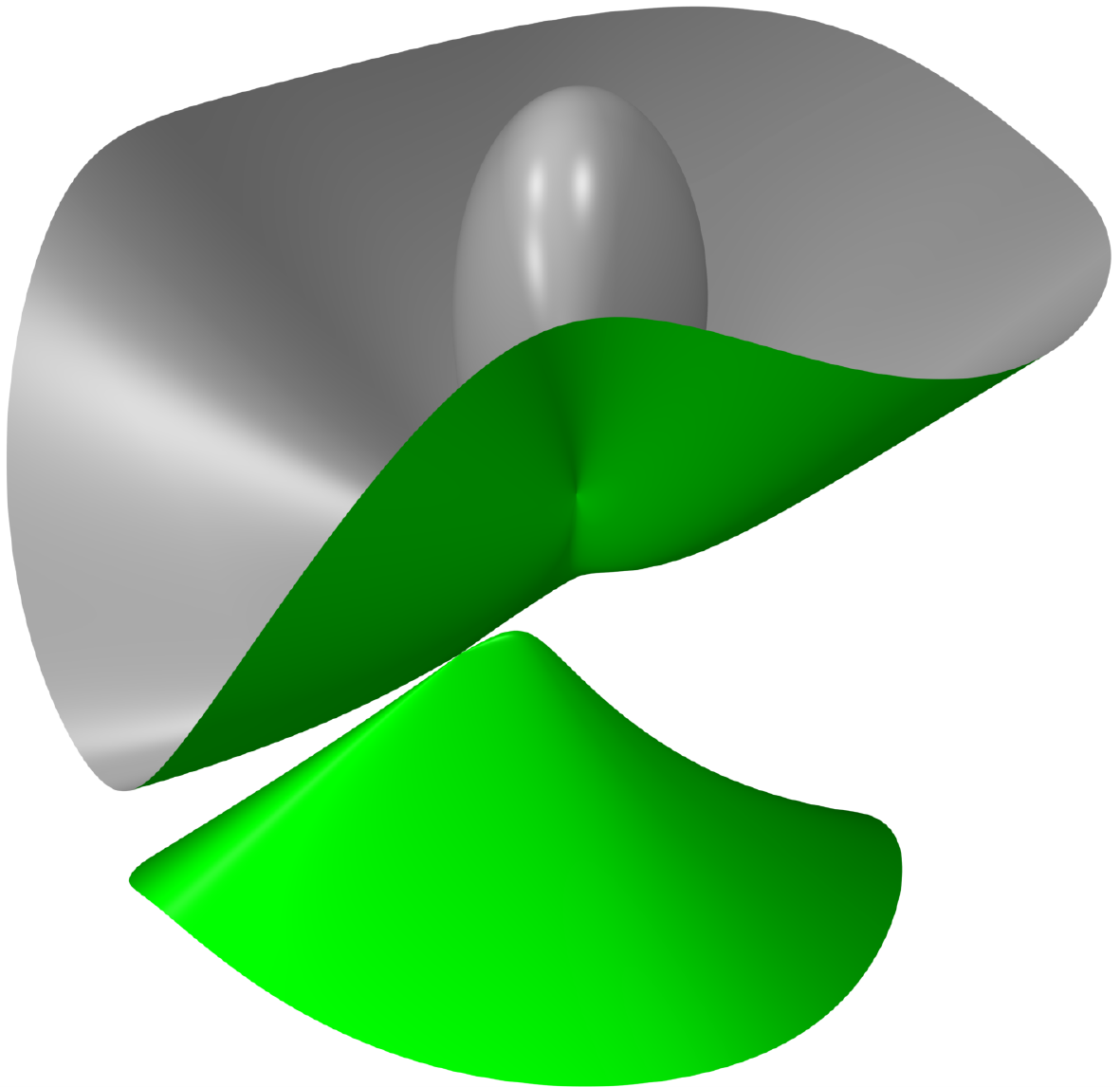}
    \caption{A real, quartic surface
    with a definite Hermitian determinantal representation
    and $(\eta, \rho, \sigma) = (2, 2, 1)$.}
\end{figure}

\subsection{Three \pointsec{2}s}

$(3, 3, 3)$:
\begin{equation*}
    \resizebox{\textwidth}{!}{\(
        \begin{bmatrix}
            x_0 + x_1 + x_2 + x_3  & x_1 + x_2               &
            -x_1 + x_2 + \ii x_0   & \ii x_2
            \\
            x_1 + x_2              & x_0 + 2x_1 + 2x_2 + x_3 &
            x_0 - x_1 + x_2        & x_0 + x_1 + x_2
            \\
            -x_1 + x_2 - \ii x_0   & x_0 - x_1 + x_2         &
            2x_0 + x_1 + x_2 + x_3 & x_0 + \ii x_2
            \\
            -\ii x_2               & x_0 + x_1 + x_2         &
            x_0 - \ii x_2          & x_0 + x_1 + 3x_2 + x_3
        \end{bmatrix}
    \)}
\end{equation*}
$\e \coloneqq [0 : 0 : 0 : 1]$

\medskip\noindent
$(3, 3, 2)$:
\begin{equation*}
    \begin{bmatrix}
        x_0 + x_2 + x_3 & x_0               &
        x_0             & 0
        \\
        x_0             & 2x_0 + 2x_2 + x_3 &
        x_0             & x_0 - \ii x_1
        \\
        x_0             & x_0               &
        x_0 + x_3       & x_1
        \\
        0               & x_0 + \ii x_1     &
        x_1             & x_0 + x_3
    \end{bmatrix}
\end{equation*}
$\e \coloneqq [0 : 0 : 0 : 1]$

\medskip\noindent
$(3, 3, 1)$:
\begin{equation*}
    \begin{bmatrix}
        x_2 + x_3 & x_0           & \ii x_2    & 0
        \\
        x_0       & x_2 + x_3     & x_2        & x_2 - \ii x_1
        \\
        - \ii x_2 & x_2           & 2x_2 + x_3 & x_1 + x_2
        \\
        0         & x_2 + \ii x_1 & x_1 + x_2  & x_2 + x_3
    \end{bmatrix}
\end{equation*}
$\e \coloneqq [0 : 0 : 0 : 1]$

\medskip\noindent
$(3, 3, 0)$:
\begin{equation*}
    \begin{bmatrix}
        x_3 & x_0      & x_2       & 0
        \\
        x_0 & x_3      & \ii x_2   & -\ii x_1
        \\
        x_2 & -\ii x_2 & x_2 + x_3 & x_1
        \\
        0   & \ii x_1  & x_1       & x_3
    \end{bmatrix}
\end{equation*}
$\e \coloneqq [0 : 0 : 0 : 1]$

\medskip\noindent
$(3, 1, 1)$:
\begin{equation*}
    \begin{bmatrix}
        x_2 + 2x_3          & 0                   &
        x_3 + \ii x_2       & x_0 + \ii x_1
        \\
        0                   & x_2 + x_3           &
        x_0 + x_2 + \ii x_1 & x_2
        \\
        x_3 - \ii x_2       & x_0 + x_2 - \ii x_1 &
        2x_2 + x_3          & x_2
        \\
        x_0 - \ii x_1       & x_2                 &
        x_2                 & x_2 + x_3
    \end{bmatrix}
\end{equation*}
$\e \coloneqq [0 : 0 : 0 : 1]$

\medskip\noindent
$(3, 1, 0)$:
\begin{equation*}
    \begin{bmatrix}
        2x_3                & 0             & x_3           & x_0 + x_2 + \ii x_1
        \\
        0                   & x_3           & x_0 + \ii x_1 & x_2
        \\
        x_3                 & x_0 - \ii x_1 & x_3           & x_2
        \\
        x_0 + x_2 - \ii x_1 & x_2           & x_2           & x_3
    \end{bmatrix}
\end{equation*}
$\e \coloneqq [0 : 0 : 0 : 1]$

\begin{figure}[htb]
    \centering
    \includegraphics[height = 0.3\textheight]{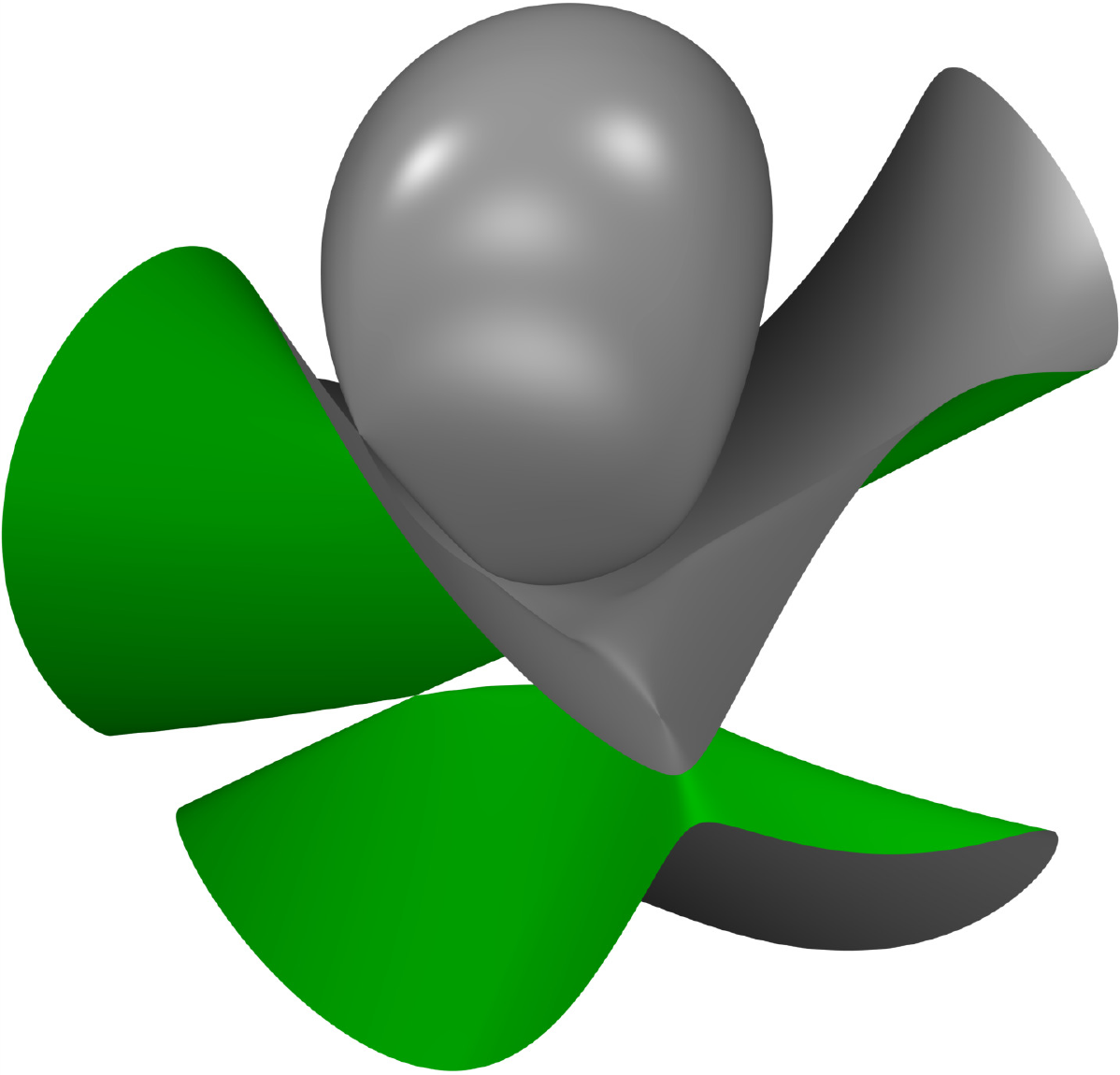}
    \caption{A real, quartic surface
    with a definite Hermitian determinantal representation
    and $(\eta, \rho, \sigma) = (3, 3, 1)$.}
\end{figure}

\subsection{Four \pointsec{2}s}

$(4, 4, 4)$:
$M_{4}(\x) \coloneqq M_{4, 0}x_0 + M_{4, 1}x_1 + M_{4, 2}x_2 + M_{4, 3}x_3$,
where
\begin{alignat*}{2}
    M_{4, 0}
    &\coloneqq
    \begin{bmatrix}
         1 &   9 &  24 &  18
        \\
         9 & -23 &   0 & -30
        \\
        24 &   0 & -36 &  36
        \\
        18 & -30 &  36 &   0
    \end{bmatrix}
    &&M_{4, 1}
    \coloneqq
    \begin{bmatrix*}[r]
         1 & 0 & -3 & 0
        \\
         0 & 4 &  0 & 6
        \\
        -3 & 0 &  9 & 0
        \\
         0 & 6 &  0 & 9
    \end{bmatrix*}
    \\
    M_{4, 2}
    &\coloneqq
    \begin{bmatrix}
        31          & -9 & -12 + 36\ii & -18
        \\
        -9          & 43 & 36          & 42
        \\
        -12 - 36\ii & 36 & 72          & 0
        \\
        -18         & 42 & 0           & 0
    \end{bmatrix}
    \,\,\,\,\,
    &&M_{4, 3}
    \coloneqq
    \begin{bmatrix}
            2 &  3 &     6 & 3\ii
        \\
            3 &  2 &     3 &   -3
        \\
            6 &  3 &    -6 & 3\ii
        \\
        -3\ii & -3 & -3\ii &    0
    \end{bmatrix}
\end{alignat*}
$\e \coloneqq [0 : 12 : 0 : 1]$

\medskip\noindent
$(4, 4, 3)$:
$M_{4}(\x) \coloneqq M_{4, 0}x_0 + M_{4, 1}x_1 + M_{4, 2}x_2 + M_{4, 3}x_3$,
where
\begin{alignat*}{2}
    M_{4, 0}
    &\coloneqq
    \begin{bmatrix}
        26 &  0 & 12 & 18
        \\
         0 & -4 &  0 & -6
        \\
        12 &  0 & -9 &  9
        \\
        18 & -6 &  9 &  0
    \end{bmatrix}
    &&M_{4, 1}
    \coloneqq
    \begin{bmatrix*}[r]
         1 & 0 & -3 & 0
        \\
         0 & 4 &  0 & 6
        \\
        -3 & 0 &  9 & 0
        \\
         0 & 6 &  0 & 9
    \end{bmatrix*}
    \\
    M_{4, 2}
    &\coloneqq
    \begin{bmatrix}
        -2       & 0 & -1 + \ii & -2
        \\
        0        & 1 & 1        & 1
        \\
        -1 - \ii & 1 & 2        & 0
        \\
        -2       & 1 & 0        & 0
    \end{bmatrix}
    \qquad
    &&M_{4, 3}
    \coloneqq
    \begin{bmatrix}
            2 &  3 &     6 & 3\ii
        \\
            3 &  2 &     3 &   -3
        \\
            6 &  3 &    -6 & 3\ii
        \\
        -3\ii & -3 & -3\ii &    0
    \end{bmatrix}
\end{alignat*}
$\e \coloneqq [0 : 4 : 0 3]$

\medskip\noindent
$(4, 4, 2)$:
\begin{equation*}
    \begin{bmatrix}
        x_0 + x_3 & x_0        & x_2       & 0
        \\
        x_0       & 2x_0 + x_3 & \ii x_2   & -\ii x_1
        \\
        x_2       & -\ii x_2   & x_2 + x_3 & x_1
        \\
        0         & \ii x_1    & x_1       & x_3
    \end{bmatrix}
\end{equation*}
$\e \coloneqq [0 : 0 : 0 : 1]$

\medskip\noindent
$(4, 4, 1)$:
\begin{equation*}
    \begin{bmatrix}
        x_0 + 2x_2 - 5x_3  & -x_0 + x_2 - 2x_3  &
        x_0 + x_3          & 2x_0 + \ii x_3
        \\
        -x_0 + x_2 - 2x_3  & 2x_1 - 2x_2 + 8x_3 &
        2x_1 - 2x_2 + 7x_3 & x_1 - x_2 + 4x_3
        \\
        x_0 + x_3          & 2x_1 - 2x_2 + 7x_3 &
        2x_1 - 2x_2 + 7x_3 & x_1 + \ii x_3
        \\
        2x_0 - \ii x_3     & x_1 - x_2 + 4x_3   &
        x_1 - \ii x_3      & x_2
    \end{bmatrix}
\end{equation*}
$\e \coloneqq [0 : 5 : 32 : 10]$

\medskip\noindent
$(4, 4, 0)$:
\begin{equation*}
    \begin{bmatrix}
        x_3           & x_0     & x_1 + \ii x_2 & 0
        \\
        x_0           & x_3     & 0             & -\ii x_1
        \\
        x_1 - \ii x_2 & 0       & x_3           & x_1 + \ii x_2
        \\
        0             & \ii x_1 & x_1 -\ii x_2  & x_3
    \end{bmatrix}
\end{equation*}
$\e \coloneqq [0 : 0 : 0 : 1]$

\medskip\noindent
$(4, 2, 2)$:
\begin{equation*}
    \begin{bmatrix}
        -77x_1 + x_2          & 27x_1                 &
        -12x_1 - 3x_2         & x_0 + 81x_1 + \ii x_3
        \\
        27x_1                 & -74x_1 + 4x_2         &
        x_0 + \ii x_3         & 78x_1 + 6x_2
        \\
        -12x_1 - 3x_2         & x_0 - \ii x_3         &
        -45x_1 + 9x_2         & 54x_1
        \\
        x_0 + 81x_1 - \ii x_3 & 78x_1 + 6x_2          &
        54x_1                 & 9x_2
    \end{bmatrix}
\end{equation*}
$\e \coloneqq [0 : -1 : 18 : 1]$

\medskip\noindent
$(4, 2, 1)$:
\begin{equation*}
    \begin{bmatrix}
        26x_1 + x_2           & 0                     &
        12x_1 - 3x_2          & x_0 + 18x_1 + \ii x_3
        \\
        0                     & -4x_1 + 4x_2          &
        x_0 + \ii x_3         & -6x_1 + 6x_2
        \\
        12x_1 - 3x_2          & x_0 - \ii x_3         &
        -9x_1 + 9x_2          & 9x_1
        \\
        x_0 + 18x_1 - \ii x_3 & -6x_1 + 6x_2          &
        9x_1                  & 9x_2
    \end{bmatrix}
\end{equation*}
$\e \coloneqq [-45 : 1 : 201 : 9]$

\medskip\noindent
$(4, 0, 0)$:
\begin{equation*}
    \begin{bmatrix}
        2x_3          & x_2           & x_2 + x_3     & x_0 + \ii x_1
        \\
        x_2           & x_3           & x_0 + \ii x_1 & 0
        \\
        x_2 + x_3     & x_0 - \ii x_1 & x_2 + x_3     & 0
        \\
        x_0 - \ii x_1 & 0             & 0             & x_3
    \end{bmatrix}
\end{equation*}
$\e \coloneqq [0 : 0 : 0 : 1]$

\begin{figure}[hbt]
    \centering
    \includegraphics[height = 0.22\textheight]{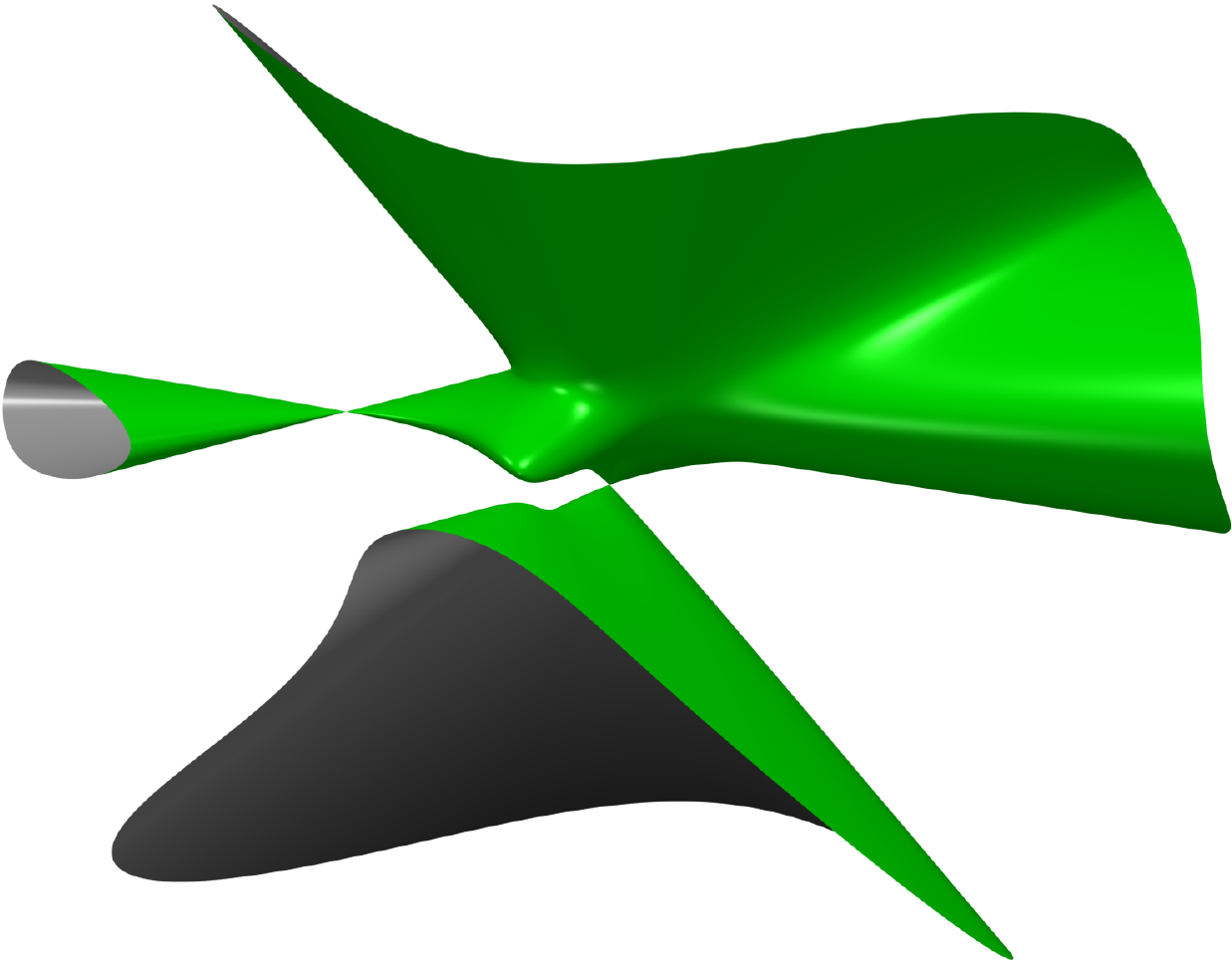}%
    \hfill
    \includegraphics[height = 0.22\textheight]{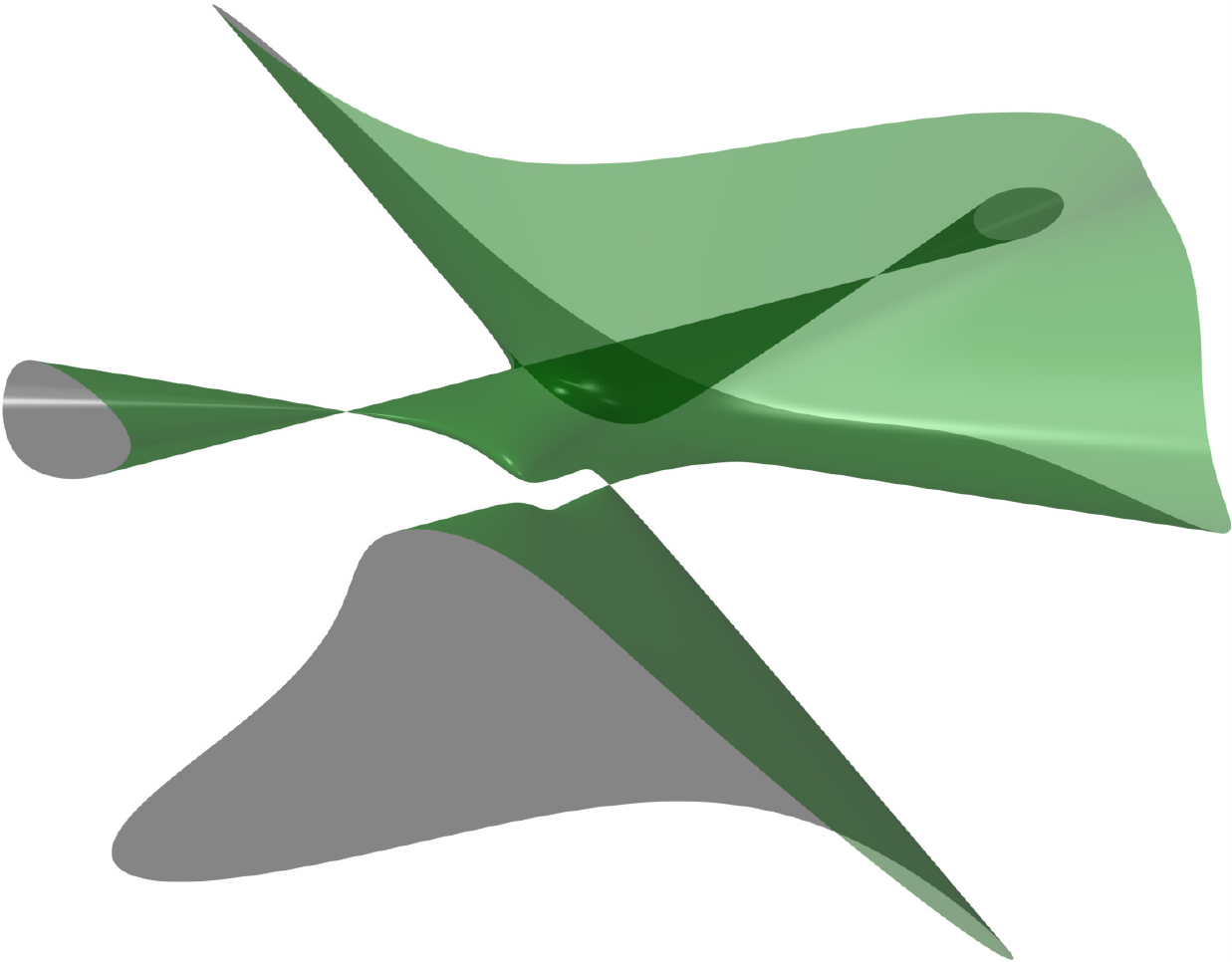}
    \caption{A real, quartic surface
        with a definite Hermitian determinantal representation
        and $(\eta, \rho, \sigma) = (4, 4, 2)$.}
\end{figure}

\subsection{Five \pointsec{2}s}

$(5, 5, 3)$:
\begin{equation*}
    \begin{bmatrix}
        x_0 + x_3 & \ii x_0   & x_2       & 0
        \\
        -\ii x_0  & x_0 + x_3 & \ii x_2   & -\ii x_1
        \\
        x_2       & -\ii x_2  & x_2 + x_3 & x_1
        \\
        0         & \ii x_1   & x_1       & x_3
    \end{bmatrix}
\end{equation*}
$\e \coloneqq [0 : 0 : 0 : 1]$

\medskip\noindent
$(5, 5, 2)$:
\begin{equation*}
    \begin{bmatrix}
        x_0 + x_2 + x_3 & x_0 + x_2           &
        x_0 - x_2       & 0
        \\
        x_0 + x_2       & 2x_0 + 2x_2 + x_3   &
        x_0 - x_2       & x_0 + x_2 - \ii x_1
        \\
        x_0 - x_2       & x_0 - x_2           &
        x_0 + x_2 + x_3 & x_1
        \\
        0               & x_0 + x_2 + \ii x_1 &
        x_1             & x_0 + x_2 + x_3
    \end{bmatrix}
\end{equation*}
$\e \coloneqq [0 : 0 : 0 : 1]$

\medskip\noindent
$(5, 5, 1)$:
\begin{equation*}
    \begin{bmatrix}
        x_0 + 9x_1 - x_2 - x_3  & -x_0 + 3x_1 + x_2 - x_3 &
        x_0 + \ii x_3           & 2x_0
        \\
        -x_0 + 3x_1 + x_2 - x_3 & x_3                     &
        x_3                     & x_3
        \\
        x_0 - \ii x_3           & x_3                     &
        x_1 - x_2 + 2x_3        & 4x_1
        \\
        2x_0                    & x_3                     &
        4x_1                    & 4x_2
    \end{bmatrix}
\end{equation*}
$\e \coloneqq [0 : 21 : 22 : 80]$

\begin{figure}[hbtp]
    \centering
    \includegraphics[height = 0.27\textheight]{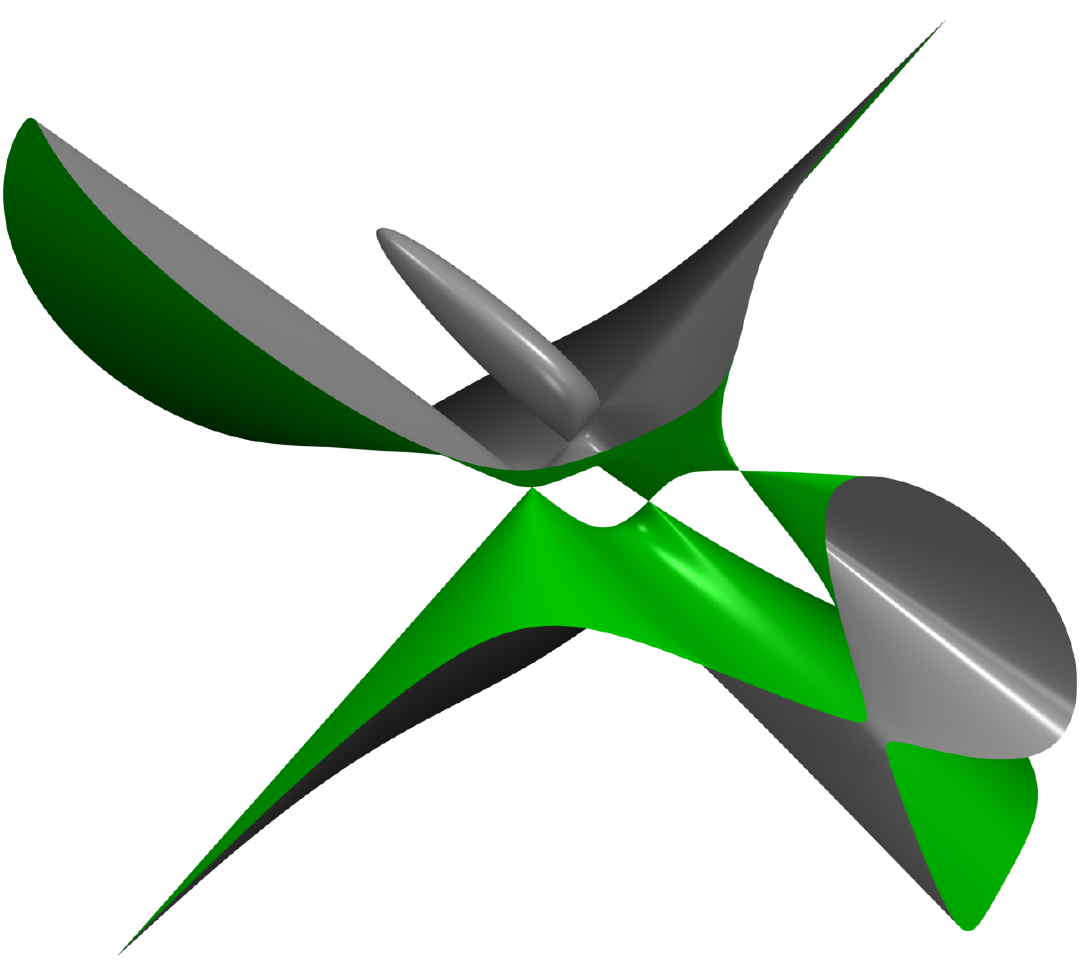}
    \caption{A real, quartic surface
    with a definite Hermitian determinantal representation
    and $(\eta, \rho, \sigma) = (5, 5, 2)$.}
\end{figure}

\subsection{Six \pointsec{2}s}

$(6, 6, 4)$:
\begin{equation*}
    \begin{bmatrix}
        x_3 & x_0     & x_1           & 0
        \\
        x_0 & x_3     & 0             & -\ii x_1
        \\
        x_1 & 0       & x_2 + x_3     & x_1 + \ii x_2
        \\
        0   & \ii x_1 & x_1 - \ii x_2 & x_3
    \end{bmatrix}
\end{equation*}
$\e \coloneqq [0 : 0 : 0 : 1]$

\medskip\noindent
$(6, 6, 3)$:
\begin{equation*}
    \begin{bmatrix}
        -x_2 + x_3 & x_0             & x_1        & \ii x_2
        \\
        x_0        & x_0 - x_2 + x_3 & x_1        & x_0
        \\
        x_1        & x_1             & -x_2 + x_3 & \ii x_2
        \\
        -\ii x_2   & x_0             & - \ii x_2  & x_3
\end{bmatrix}
\end{equation*}
$\e \coloneqq [0 : 0 : 1 : 3]$

\begin{figure}[htbp]
    \centering
    \includegraphics[height = 0.25\textheight]{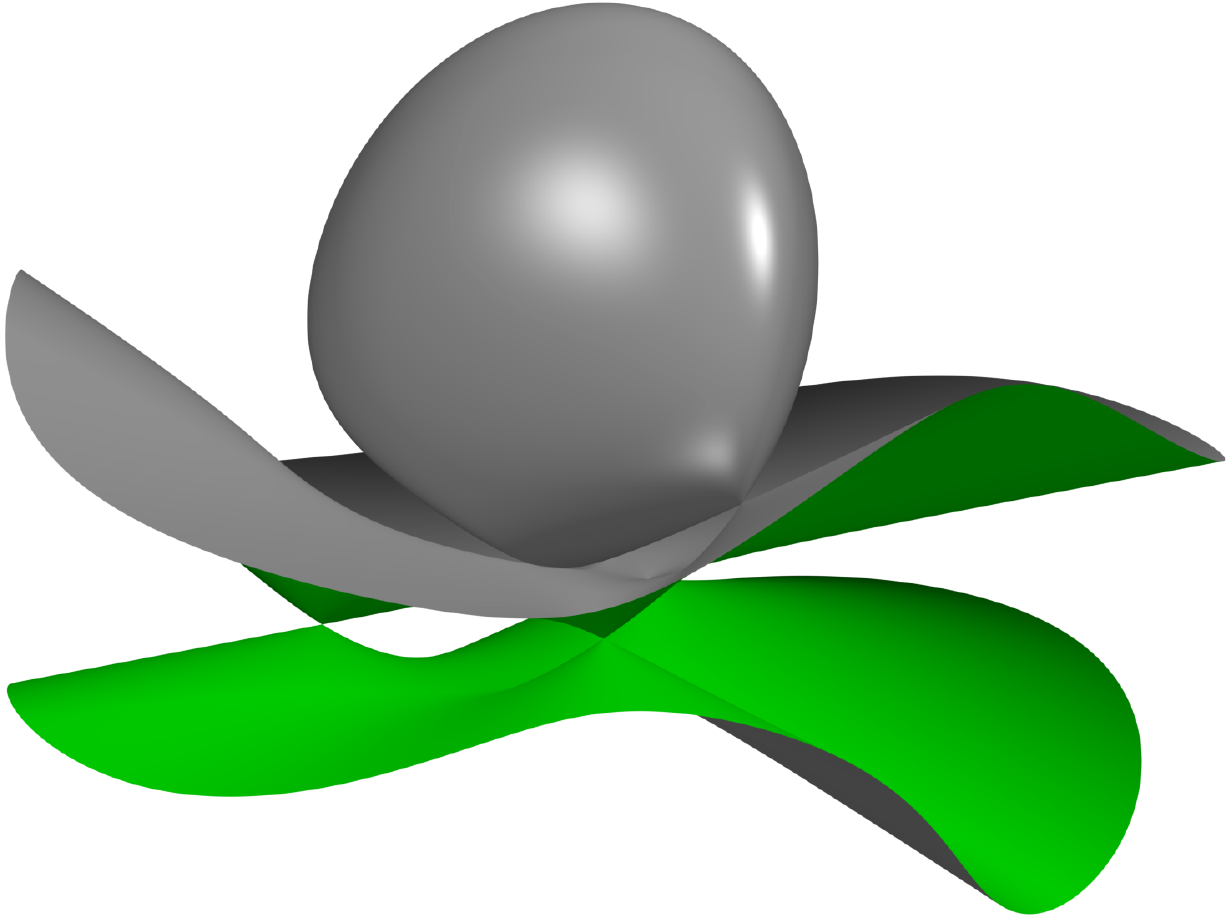}
    \caption{A real, quartic surface
    with a definite Hermitian determinantal representation
    and $(\eta, \rho, \sigma) = (6, 6, 3)$.
    Not all nodes are visible from this angle.}
\end{figure}

\subsection{Seven \pointsec{2}s}

$(7, 7, 5)$:
\begin{equation*}
    \begin{bmatrix}
        x_0 + x_3 & 0         & x_2     & 0
        \\
        0         & x_0 + x_3 & \ii x_2 & -\ii x_1
        \\
        x_2       & -\ii x_2  & x_2+x_3 & x_1
        \\
        0         & \ii x_1   & x_1     & x_3
    \end{bmatrix}
\end{equation*}
$\e \coloneqq [0 : 0 : 0 : 1]$

\medskip\noindent
$(7, 7, 4)$:
\begin{equation*}
    \begin{bmatrix}
        -x_2 + x_3                 & x_0 - \ii(x_0 i + x_2 - x_3) &
        x_1                        & -x_0 + x_3
        \\
        x_0 + \ii(x_0 + x_2 - x_3) & -x_0 - x_2 + 2x_3            &
        x_1                        & x_0 + \ii(x_0 + x_2 - x_3)
        \\
        x_1                        & x_1                          &
        x_0                        & x_2
        \\
        -x_0 + x_3                 & x_0 - \ii(x_0 + x_2 - x_3)   &
        x_2                        & x_3
    \end{bmatrix}
\end{equation*}
$\e \coloneqq [1 : 0 : 0 : 2]$

\begin{figure}[hbtp]
    \centering
    \includegraphics[height = 0.25\textheight]{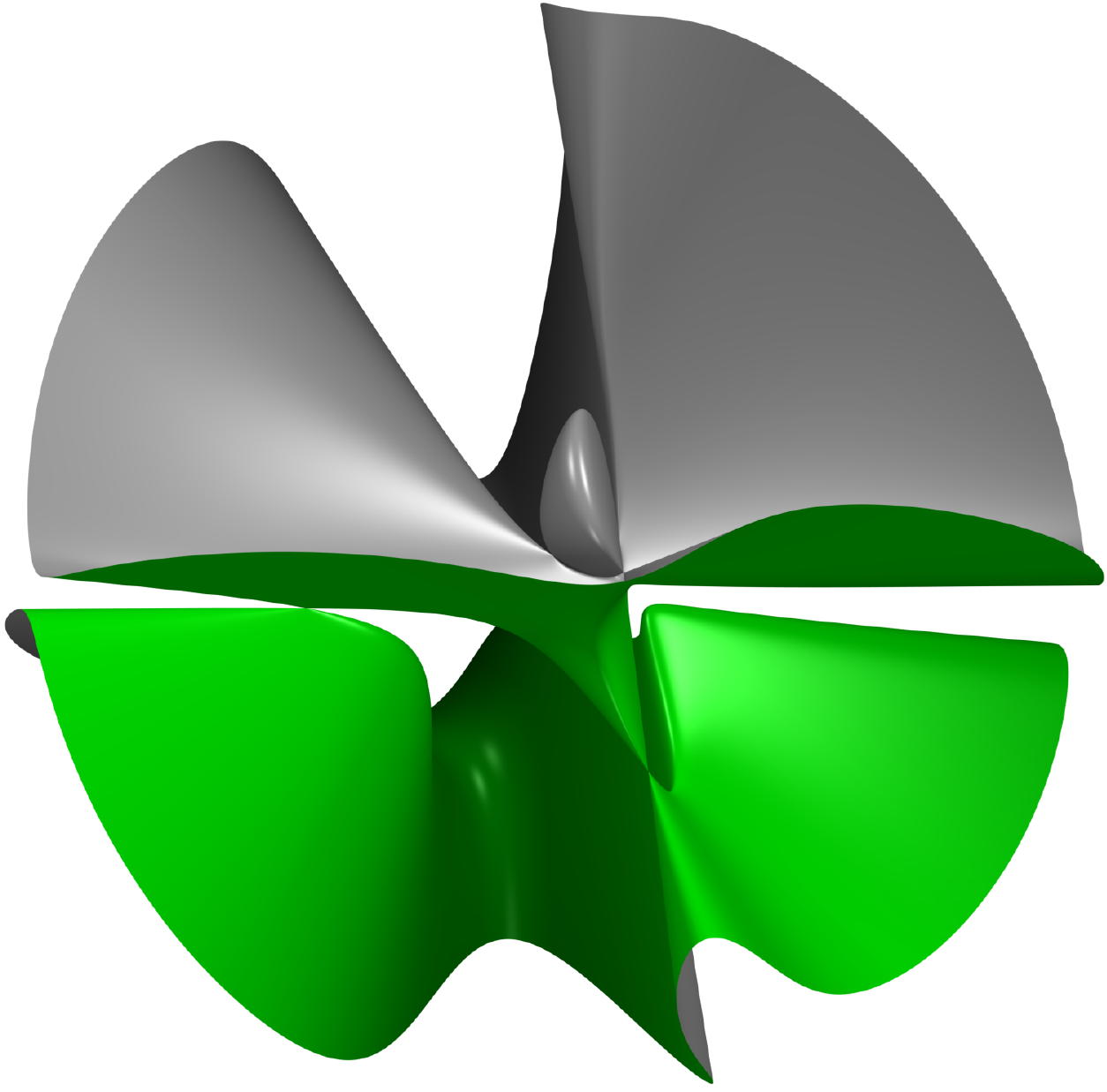}
    \addtolength{\belowcaptionskip}{-2em}
    \caption{A real, quartic surface
    with a definite Hermitian determinantal representation
    and $(\eta, \rho, \sigma) = (7, 7, 4)$.
    Not all nodes are visible from this angle.}
\end{figure}

\subsection{Eight \pointsec{2}s}

$(8, 8, 5)$:
\begin{equation*}
    \begin{bmatrix}
        -x_2 + x_3                 & x_0 - \ii(x_0 + x_2 - x_3) &
        x_1                        & -x_0 + x_3
        \\
        x_0 + \ii(x_0 + x_2 - x_3) & 2x_0                       &
        x_1                        & x_0 + \ii(x_0 + x_2 - x_3)
        \\
        x_1                        & x_1                        &
        x_0                        & x_2
        \\
        -x_0 + x_3                 & x_0 - \ii(x_0 + x_2 - x_3) &
        x_2                        & x_3
    \end{bmatrix}
\end{equation*}
$\e \coloneqq [1 : 0 : 0 : 2]$

\medskip\noindent
$(8, 8, 4)$:
\begin{equation*}
    \begin{bmatrix}
        2x_0 + x_3 & \ii x_0   & x_2       & 0
        \\
        -\ii x_0   & x_0 + x_3 & \ii x_2   & -\ii x_1
        \\
        x_2        & -\ii x_2  & x_2 + x_3 & x_1
        \\
        0          & \ii x_1   & x_1       & x_3
    \end{bmatrix}
\end{equation*}
$\e \coloneqq [0 : 0 : 0 : 1]$

\medskip\noindent
$(8, 6, 4)$:
\begin{equation*}
    \begin{bmatrix}
        x_0 + x_3 & 0          & x_2       & 0
        \\
        0         & 2x_0 + x_3 & \ii x_2   & -\ii x_1
        \\
        x_2       & -\ii x_2   & x_2 + x_3 & x_1
        \\
        0         & \ii x_1    & x_1       & x_3
    \end{bmatrix}
\end{equation*}
$\e \coloneqq [0 : 0 : 0 : 1]$

\begin{figure}[htbp]
    \centering
    \includegraphics[height = 0.22\textheight]{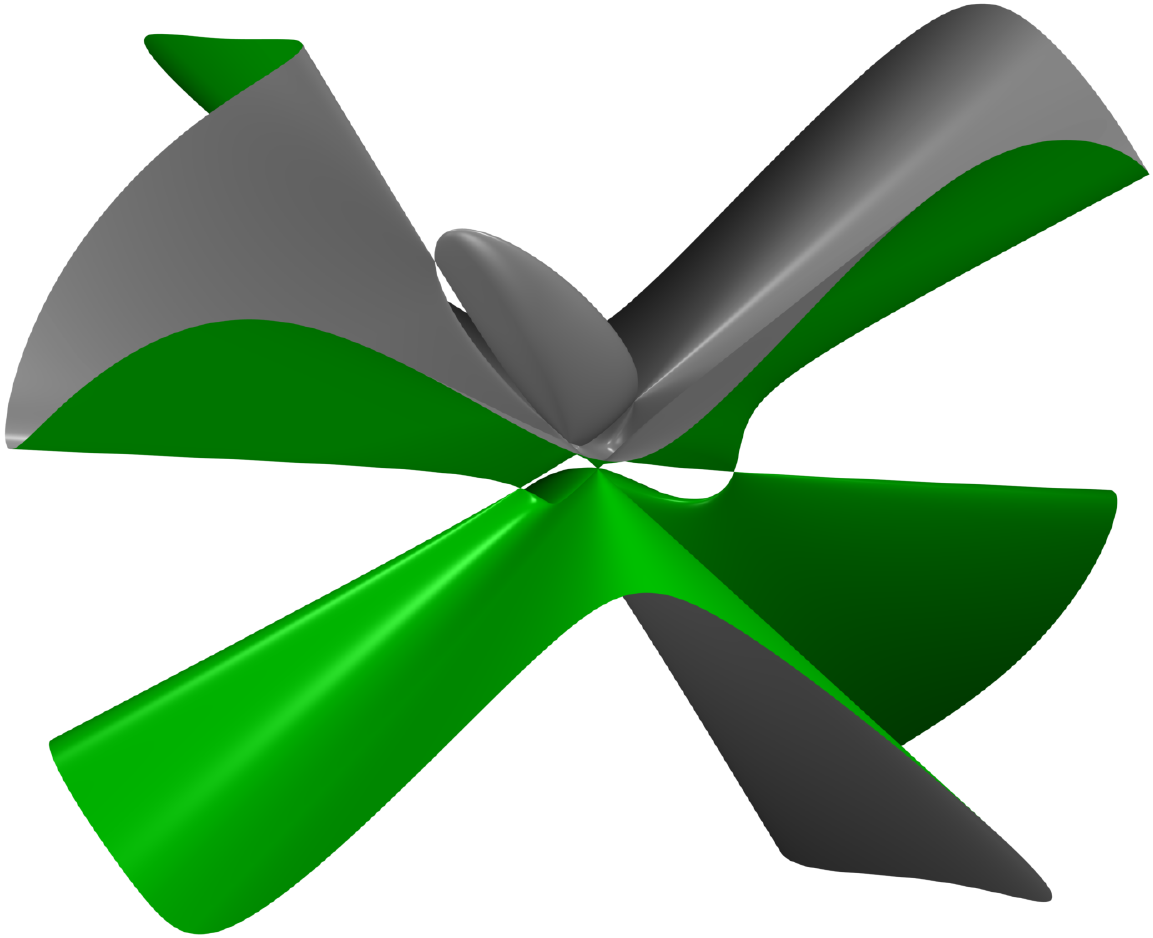}
    \caption{A real, quartic surface
    with a definite Hermitian determinantal representation
    and $(\eta, \rho, \sigma) = (8, 8, 5)$.
    Not all nodes are visible from this angle.}
\end{figure}

\printbibliography

\paragraph{Author's address:}

Martin Hels\o,
University of Oslo,
Postboks 1053 Blindern,
0316 Oslo,
Norway,
\href{mailto:martibhe@math.uio.no}{\nolinkurl{martibhe@math.uio.no}}

\end{document}